\documentclass[11pt, a4paper]{article}

\usepackage{fullpage, comment}
\usepackage{psfrag}
\usepackage{latexsym}
\usepackage{import, pst-plot}
\usepackage{amsfonts}
\usepackage{amsmath}
\usepackage{amssymb}
\usepackage{graphicx}

\newcommand{\beq}{\begin{equation}}
\newcommand{\eeq}{\end{equation}}
\newcommand{\ba}{\begin{array}}
\newcommand{\ea}{\end{array}}
\newcommand{\bea}{\begin{eqnarray}}
\newcommand{\eea}{\end{eqnarray}}
\newcommand{\bc}{\begin{center}}
\newcommand{\ec}{\end{center}}
\newcommand{\bt}{\begin{table}}
\newcommand{\et}{\end{table}}
\newcommand{\eps}{\epsilon}
\newcommand{\la}[1]{\label{#1}}
\newcommand{\p}{\partial}

\newcommand{\ds}{\displaystyle}
\newcommand{\no}{\noindent}
\newcommand{\pp}[2]{{\partial #1 \over \partial #2}}

\newcommand{\rf}[1]{(\ref{#1})}
\newcommand{\beqno}{\begin{displaymath}}
\newcommand{\eeqno}{\end{displaymath}}
\newcommand{\dd}[2]{{\delta #1 \over \delta #2}}

\newcommand{\been}{\begin{enumerate}}
\newcommand{\een}{\end{enumerate}}

\renewcommand{\iff}{\Leftrightarrow}

\newcommand{\sn}{{\rm sn}}
\newcommand{\cn}{{\rm cn}}
\newcommand{\dn}{{\rm dn}}

\newcommand{\vs}{\vspace*{0.1in}}
\newcommand{\ra}{\rightarrow}
\newcommand{\sgn}{\mbox{sign}}

\def\Xint#1{\mathchoice
   {\XXint\displaystyle\textstyle{#1}}%
   {\XXint\textstyle\scriptstyle{#1}}%
   {\XXint\scriptstyle\scriptscriptstyle{#1}}%
   {\XXint\scriptscriptstyle\scriptscriptstyle{#1}}%
   \!\int}
\def\XXint#1#2#3{{\setbox0=\hbox{$#1{#2#3}{\int}$}
     \vcenter{\hbox{$#2#3$}}\kern-.5\wd0}}

\def\dashint{\Xint-}

\newlength{\myheight}
\newlength{\mylength}

\newcounter{saveeqn}

\newcommand{\alpheqn}{\setcounter{saveeqn}{\value{equation}}
\stepcounter{saveeqn}\setcounter{equation}{0}
\renewcommand{\theequation}{\mbox{\arabic{saveeqn}\alph{equation}}}}

\newcommand{\resetalpheqn}{\setcounter{equation}{\value{saveeqn}}
\renewcommand{\theequation}{\arabic{equation}}}

\newtheorem{example}{Example}

\begin{document}

\title{High-frequency instabilities of small-amplitude solutions of Hamiltonian PDEs}

\author{
Bernard Deconinck$^\dagger$ \& Olga Trichtchenko$^*$\\
~\\
$^\dagger$ Department of Applied Mathematics\\
University of Washington\\
Campus Box 352420, Seattle, WA, 98195, USA\\
e-mail: {\tt deconinc@uw.edu}\\
~\\
$^*$ Department of Mathematics\\
University College London\\
Gower Street\\
London, WC1E 6BT, UK\\
e-mail: {\tt olga.trichtchenko@gmail.com}}

\maketitle

\begin{abstract}

Generalizing ideas of MacKay, and MacKay and Saffman, a necessary
condition for the presence of high-frequency ({\em i.e.}, not modulational)
instabilities of small-amplitude periodic solutions of Hamiltonian partial
differential equations is presented, entirely in terms of the
Hamiltonian of the linearized problem. With the
exception of a Krein signature calculation, the theory is completely phrased in terms of
the dispersion relation of the linear problem. The general theory
changes as the Poisson structure of the Hamiltonian partial differential equation is changed.
Two important cases of such Poisson structures are worked out in full generality. An example not fitting these two important cases is
presented as well, using a candidate Boussinesq-Whitham equation.

\end{abstract}

\section{Introduction}

It is expected that much of the dynamics of the small-amplitude solutions of a partial differential equation (PDE), including their stability or instability, is dictated by the study of a linearized (about a trivial solution, say $u=0$) problem. In this article, we focus specifically on the spectral stability of periodic traveling-wave solutions of Hamiltonian PDEs as they bifurcate away from a trivial solution. Our work follows earlier ideas of MacKay \cite{mackay} and MacKay and Saffman \cite{mackaysaffman}. We start from an autonomous Hamiltonian system of PDEs \cite{arnold}, {\em i.e.},

\beq\la{ham}
u_t=J \dd{H}{u}.
\eeq

\no \sloppypar \no Here and throughout, indices involving $x$ or $t$ denote partial derivatives. Further, $u=(u_1(x,t), \ldots, u_M(x,t))^T$ is an $M$-dimensional vector function defined in a suitable function space, and $J$ is a Poisson operator \cite{arnold, arnoldnovikov}. More details and examples are given below. Finally, $H=\int_D {\cal H}(u,u_x, \ldots)dx$ is the Hamiltonian, whose density ${\cal H}$ depends on $u$ and its spatial derivatives, defined for $x\in D$. We consider only the stability of periodic solutions, thus $D$ is any interval of length $L$, the period.
Note that for some of our examples ${\cal H}$ will depend on spatial derivatives of $u$ of arbitrary order.

To investigate the stability of traveling wave solutions of this system, we reformulate \rf{ham} in a frame moving with speed $c$, using the transformation $\hat x=x-ct$, $\hat t=t$, and considering solutions $u(\hat x, \hat t)=U(\hat x)$ (successively omitting hats). This leads to

\beq\la{hamc}
u_t-c u_x=J \dd{H}{u} ~~\iff~~u_t=J \dd{H_c}{u},
\eeq

\no for a modified Hamiltonian $H_c$. Traveling wave solutions are solutions of the ordinary differential system

\beq\la{travham}
-c U_x=J \dd{H}{U}  ~~\iff~~0=J \dd{H_c}{U}.
\eeq

\no Thus if $J$ is invertible, traveling waves are stationary points of the Hamiltonian $H_c$. The system~\rf{travham} typically has the zero (trivial) solution for a range of $c$ values. The small-amplitude solutions whose stability we investigate bifurcate away from these trivial zero-amplitude solutions at special values of the speed parameter $c$, as is schematically shown in Fig.~\ref{fig:bif}. It is our goal to see to what extent anything can be said about the stability of the small-amplitude solutions (with amplitudes in the shaded regions of Fig.~\ref{fig:bif}) from knowledge of the zero-amplitude solutions at the bifurcation point. An outline of the steps in this process is as follows.

%
\begin{figure}[tb]
\def\svgwidth{5in}
\centerline{\hspace*{0.3in}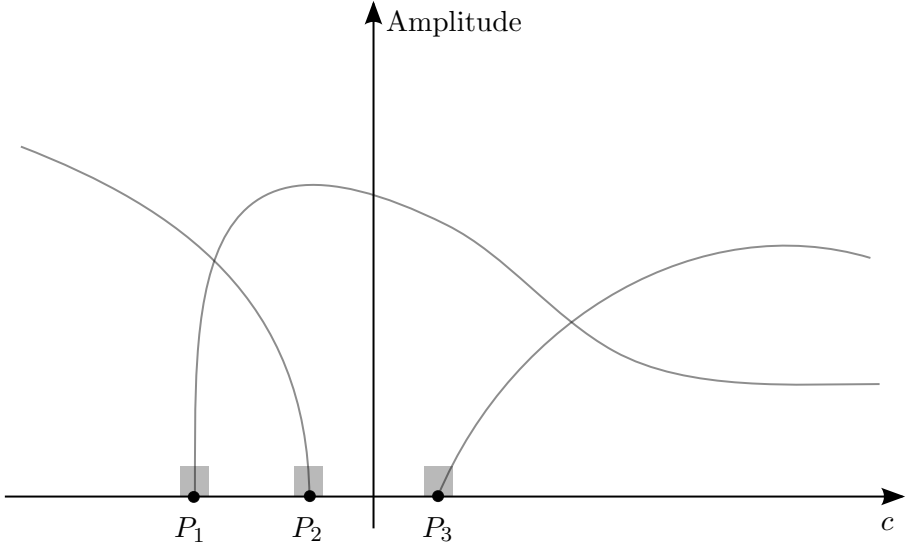}
\caption{\la{fig:bif} A cartoon of the bifurcation structure of the traveling waves for a third-order ($M=3$) system: solution branches bifurcate away from the trivial zero-amplitude solution at specific values of the traveling wave speed $c$.
    }
\end{figure}
%

\begin{enumerate}

\item {\bf Quadratic Hamiltonian}. A linear system of equations is obtained by linearizing the system \rf{travham} around the zero solution: let $u=\epsilon v+{o}(\epsilon)$ and omit terms of order $o(\epsilon)$. Alternatively, if $J$ is independent of $u$ and its spatial derivatives, one may expand the Hamiltonian $H_c$ as a function of $\epsilon$ and retain its quadratic terms. The resulting Hamiltonian $H_c^0$ of the linearized system is the starting point for the next steps.

\item {\bf Dispersion relation}. The linearized system has constant coefficients and is easily solved using Fourier analysis. The dispersion relation $F(\omega, k)=0$ governs the time dependence of the solutions. It is obtained by investigating solutions whose spatial and temporal dependence is proportional to $\exp(ikx-i \omega t)$. Here $F(\omega,k)=0$ is of degree $M$ in $\omega$. It is a fundamental assumption of our approach that all solutions $\omega_j(k)$ ($j=1, \ldots, M$) of $F(\omega,k)=0$ are real for $k\in \mathbb{R}$. The dispersion relation can be expressed entirely in terms of the coefficients appearing in the quadratic Hamiltonian $H_c^0$. For periodic systems of period $L$, the values of $k$ are restricted to be of the form $2\pi N/L$, $N\in \mathbb{Z}$.

\item {\bf Bifurcation branches}. The values of the phase speed $c_j=\omega_j/k$ for which nontrivial solutions bifurcate away from the zero-amplitude solution are determined by the condition that the zero solution is {\em not} the unique solution to the Fourier transformed problem. In effect, this is the classical bifurcation condition that a Jacobian is singular. This simple calculation determines the bifurcation branch starting points explicitly in terms of the different solutions to the dispersion relation. In what follows, we follow the first branch, starting at $c_1$, without loss of generality.

    It is assumed that only a single non-trivial bifurcation branch emanates from a bifurcation point. Although more general cases can be incorporated, we do not consider them here. Further, we fix the period of the solutions on the bifurcation branch (usually to $2\pi$). Other choices can be made. Instead of varying the amplitude as a function of the speed for fixed period, one could fix the speed and vary the period, {\em etc.} The methods presented can be redone for those scenarios in a straightforward fashion.

\item {\bf Stability spectrum}. The spectrum of the linear operator determining the spectral stability of the zero solution at the bifurcation point on the first branch is calculated. Since this spectral problem has constant coefficients, this calculation can be done explicitly. Again, this is done entirely in terms of the dispersion relation of the problem. Using a Floquet decomposition (see \cite{deconinckkutz1, kapitulapromislow}), the spectrum is obtained as a collection of point spectra, parameterized by the Floquet exponent $\mu\in (-\pi/L, \pi/L]$. Due to the reality of the branches of the dispersion relation, the spectrum is confined to the imaginary axis. In other words, the zero-amplitude solutions are spectrally stable. The use of the Floquet decomposition allows for the inclusion of perturbations that are not necessarily periodic with period $L$. Instead, the perturbations may be quasiperiodic with two incommensurate periods, subharmonic (periodic, but with period an integer multiple of $L$), or spatially localized \cite{deconinckkutz1, haraguskapitula, kapitulapromislow}.

\item {\bf Collision condition}. Given the explicit expression for individual eigenvalues $\lambda$, it is easy to find the conditions for which eigenvalues corresponding to different parameters (Floquet exponent, branch number of the dispersion relation, {\em etc}.) coincide on the imaginary axis. This is referred to as the {collision condition}. Once again, it is given entirely in terms of the dispersion relation.

    It is a consequence of the Floquet theorem \cite{coddingtonlevinson} that collisions need to be considered only for spectral elements corresponding to the same value of the Floquet exponent since the subspaces of eigenfunctions for a fixed Floquet exponent are invariant under the flow of the linearized equation.

\item {\bf Krein signature}. Having obtained the stability spectrum at the starting point of the bifurcation branches, we wish to know how the spectrum evolves as we move up a bifurcation branch. One tool to investigate this is the Krein signature \cite{kollarmiller, krein1, krein2, mackay, meiss}. In essence, the Krein signature of an eigenvalue is the sign of the Hamiltonian of the linearized system evaluated on the eigenspace of the eigenvalue. Different characterizations are given below. If two imaginary eigenvalues of the same signature collide as a parameter changes, their collision does not result in them leaving the imaginary axis. Thus the collision of such eigenvalues does not result in the creation of unstable modes. In other words, it is a necessary condition for collisions to lead to instability that the Krein signature of the colliding eigenvalues is different. This scenario is illustrated in Fig.~\ref{fig:krein}. That figure also illustrates the quadrifold symmetry of the stability spectrum of the solution of a Hamiltonian system: for each eigenvalue $\lambda\in \mathbb{C}$, $\lambda^*$, $-\lambda$ and $-\lambda^*$ are also eigenvalues. Here $\lambda^*$ denotes the complex conjugate of $\lambda$. It should be noted that the occurrence of a collision is required for eigenvalues to leave the imaginary axis, due to the quadrifold symmetry of the spectrum.

%
\begin{figure}[tb]
\def\svgwidth{6in}
\centerline{\hspace*{1in}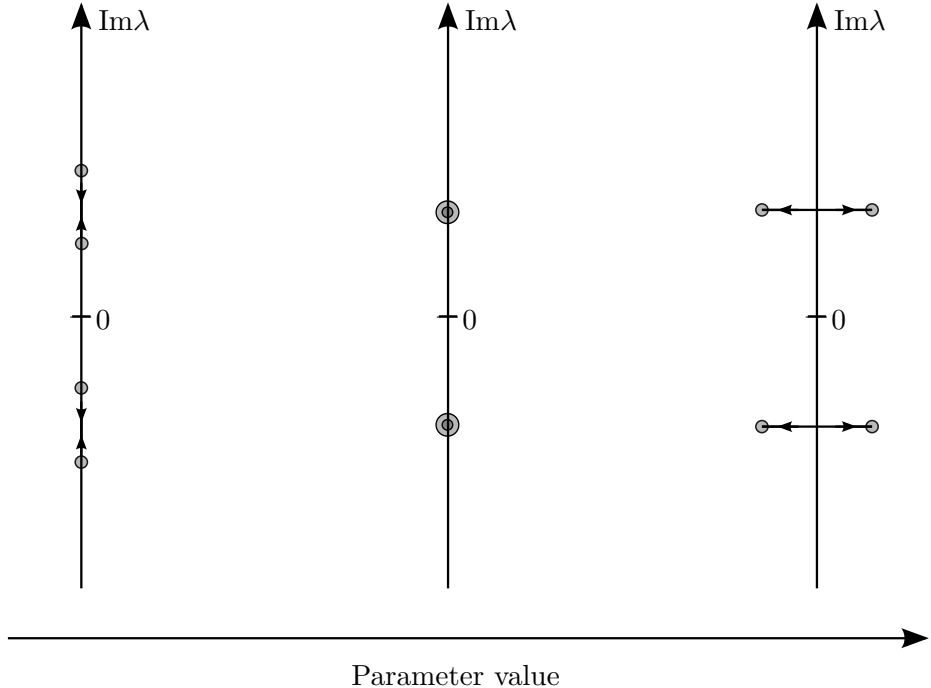}
\caption{\la{fig:krein} Colliding eigenvalues in the complex plane as a parameter is increased. On the left, two eigenvalues are moving towards each other on the positive imaginary axis, accompanied by a complex conjugate pair on the negative imaginary axis. In the middle, the eigenvalues in each pair have collided. On the right, a Hamiltonian Hopf bifurcation occurs: the collided eigenvalues separate, leaving the imaginary axis (implying that the two Krein signatures were different).
    }
\end{figure}
%

    Thus we calculate the Krein signature of any coinciding eigenvalues, obtained in Step~5. If these Krein signatures are equal, the eigenvalues will remain on the imaginary axis as the amplitude is increased. Otherwise, the eigenvalues may leave the imaginary axis, through a so-called Hamiltonian Hopf bifurcation \cite{vandermeer}, resulting in instability. Thus we establish a necessary condition for the instability of periodic solutions of small amplitude. The Krein signature condition cannot be expressed entirely in terms of the dispersion relation, and the coefficients of $H^0_c$ are required as well. Please refer to the next two sections for details.

\end{enumerate}

Although all calculations are done for the zero-amplitude solutions at the starting point of a bifurcation branch, the continuous dependence of the stability spectrum on the parameters in the problem \cite{hislop}, including the velocity of the traveling wave or the amplitude of the solutions, guarantees that the stability conclusions obtained persist for solutions of small amplitude. Thus meaningful conclusions about solutions in the shaded regions of the bifurcation branches of Fig.~\ref{fig:bif} are reached, despite the fact that we are unable to say anything about the size of the maximal amplitude for which these conclusions are valid.

\vs

{\bf Remarks.}

\begin{itemize}

\item Throughout this introduction, and in fact throughout much of this manuscript, our emphasis is on generality and usability as opposed to rigor. As a consequence some of the statements made above are necessarily vague: more precise statements would limit the generality aimed for. Within the context of more specific examples, more precision may be possible. Along the same lines, we have limited ourselves to the most generic case of two eigenvalues colliding on the imaginary axis. Lastly, since we are assuming an initial ({\em i.e.}, for a starting value of 0 for the amplitude parameter) situation that is neutrally stable, our only interest is in collisions on the imaginary axis. If eigenvalues are present off the imaginary axis for the zero-amplitude solution, such solutions are already spectrally unstable, and the continuous dependence of the spectrum on its parameters \cite{hislop} guarantees instability for solutions of small amplitude as  a consequence.

\item If eigenvalues collide at the origin, their Krein signature is zero (this follows immediately from its definition) and no conclusion can be drawn about whether or not colliding eigenvalues leave the imaginary axis or not. Because of this, the methods discussed here allow only for the study of so-called { high-frequency instabilities.} Indeed, if eigenvalues collide away from the origin, the unstable perturbations do not only grow in magnitude exponentially, they also display an oscillation of frequency $\omega$ that is equal to the non-zero imaginary part of the colliding eigenvalues. In this sense {high frequency} is more accurately described as non-zero frequency. We use the high-frequency name to distinguish the instabilities investigated here from the modulational instability \cite{zakharovostrovsky}, which is a consequence of eigenvalues colliding at the origin. Further, as seen in the examples below, often the presence of one high-frequency instability is accompanied by a sequence of such instabilities with increasing frequencies tending to infinity.

    The study of the modulational instability requires a different set of techniques, see for instance \cite{benjaminperiodic, whitham1, zakharovostrovsky} for different classical approaches. More recently, a general framework was developed by Bronski, Johnson and others (see {\em e.g.} \cite{bronskijohnson1, hurjohnson}). Collisions at the origin are common, since any Lie symmetry of the underlying problem gives rise to a zero eigenvalue in the stability problem \cite{gss1, kapitulapromislow}. For this same reason, such collisions typically involve more than two eigenvalues (for instance, the original Benjamin-Feir instability involves four \cite{benjaminperiodic, bridgesmielke}), which is one of the reasons why their treatment is often complicated: its analysis is not only technical but also involves tedious calculations, as exemplified recently in a paper by Hur and Johnson \cite{hurjohnson} on the modulational instability for periodic solutions of the Whitham equation.

\item Following the Floquet decomposition, we use eigenfunctions given by a single Fourier mode since the linear stability problem for the zero-amplitude solutions has constant coefficients. Lastly, the calculation of the Krein signature involves only the finite-dimensional eigenspace of the eigenvalue under consideration. Thus, in essence, all the calculations done in this paper are finite dimensional, as they are in \cite{mackaysaffman}, for instance. As was noted there, for equations with real-valued solutions, it is necessary to consider the eigenfunctions corresponding to Floquet exponent $\mu$ and $-\mu$ simultaneously, since the eigenfunctions with $-\mu$ are the complex conjugates of those with $\mu$.

\item We call an equation dispersive if all branches of its dispersion relation $\omega(k)$ are real valued for $k\in \mathbb{R}$. It is easy to see that there exist linear, constant-coefficient Hamiltonian systems that are not dispersive. An explicit example is
    \beq
    \left\{
    \ba{rcl}
    q_t&=&q_{xx},\\
    p_t&=&-p_{xx}.
    \ea
    \right.
    \eeq
This (admittedly bizarre) system is Hamiltonian with canonical Poisson structure
\beq\la{canonicalJ}
J=\left(
\ba{rr}
0 & 1\\
-1 & 0
\ea
\right),
\eeq
and Hamiltonian $H=-\int_0^{2\pi} q_x p_x dx$. The branches of its dispersion relation are given by $\omega_{1,2}(k)=\mp i k^2$.

It is an interesting question whether there exist dispersive systems that are not Hamiltonian. One may be tempted to consider an example like
\beq
    \left\{
    \ba{rcl}
    q_t&=&a q_{x},\\
    p_t&=&b p_{x}.
    \ea
    \right.
    \eeq
For $a, b \in \mathbb{R}$ this system is dispersive ($\omega_1=-a k$, $\omega_2=-b k$), but it is not Hamiltonian with canonical structure. However, allowing for the noncanonical structure
\beq
J=\left(
\ba{rr}
\p_x & 0\\
0 & \p_x
\ea
\right),
\eeq
the system is easily found to have the Hamiltonian $H=\frac{1}{2}\int_0^{2\pi}(a q^2+b p^2)dx$. We conjecture that no systems exist that are dispersive but not Hamiltonian. As demonstrated by the example above, the question is difficult to analyze, since different forms of the Poisson operator have to be allowed.

\item The investigation of the Hamiltonian structure of the zero-amplitude solution follows earlier work by Zakharov, Kuznetsov, and others. A review is available in \cite{zakharovmusher}. The stability of the trivial solution is also investigated there, with stability criteria given entirely in terms of the branches of the dispersion relation. Using canonical perturbation theory, Hamiltonians containing higher-than-quadratic terms are considered. This is used to consider the impact of nonlinear effects on the dynamics of the trivial solution. At this point, connections to resonant interaction theory \cite{benneyRIT, hammackhendersonRIT, phillipsRIT} become apparent, as they will in what follows. Although the physical reasoning leading to resonant interaction theory and that leading to the criteria established below is different, it is clear that connections exist. 

\end{itemize}

\section{Motivating example}\la{sec:whitham}

Our investigations began with the study of the so-called Whitham equation \cite{WhithamEq}, \cite[page 368]{whitham}. The equation is usually posed on the whole line, for which the equation satisfied by $u(x,t)$ is

\beq\la{whitham}
u_t+N(u)+\int_{-\infty}^\infty K(x-y)u_y(y,t) dy=0.
\eeq

\no The term $N(u)$ denotes the collection of all nonlinear terms in the equation. It is assumed that $\lim_{\eps\ra 0}N(\epsilon u)/\epsilon=0$. The last term encodes the dispersion relation of the linearized Whitham equation. The kernel $K(x)$ is the inverse Fourier transform of the phase speed $c(k)$, where $c(k)=\omega(k)/k$, with $\omega(k)$ the dispersion relation. Here $c(k)$ is assumed to be real valued and nonsingular for $k\in \mathbb{R}$. Thus

\beq
K(x)=\frac{1}{2\pi}\dashint_{-\infty}^\infty c(k)e^{ikx}dk,
\eeq

\no where $\dashint$ denotes the principal value integral. Depending on $c(k)$, this equation may have to be interpreted in a distributions sense \cite{EhrnstromKalish, stakgold}. Letting $u\sim \exp(ikx-\omega(k) t)$, it is a straightforward calculation to see that the dispersion relation of the linear Whitham equation is $\omega(k)$. In fact, the linear Whitham equation is easily seen to be a rewrite of the linear evolution equation \cite{as}

\beq\la{whithamasform}
u_t=-i \omega(-i \p_x)u,
\eeq

\no where $\omega(-i \p_x)$ is a linear operator with odd symbol $\omega(k)$: $\omega(k)=-\omega(-k)$, the dispersion relation considered. Indeed, letting $\omega(-i \p_x)$ act on

\beq
u(x,t)=\frac{1}{2\pi}\int_{-\infty}^\infty e^{ikx}\hat u(k,t)dk,
\eeq

\no and replacing $\hat u(k,t)$ by

\beq
\hat u(k,t)=\int_{-\infty}^\infty e^{-ik y} u(y,t) dy,
\eeq

\no the linear part of \rf{whitham} is obtained after one integration by parts. We restrict our considerations to odd dispersion relations, to ensure the reality of the Whitham equation.

One of Whitham's reasons for writing down the Whitham equation \cite{WhithamEq, whitham} was to describe waves in shallow water (leading to the inclusion of a KdV-type nonlinearity $N(u)\sim u u_x$) that feel the full dispersive response of the one-dimensional water wave problem (without surface tension), for which

\beq\la{wwdispersion}
\omega^2(k)=gk\tanh(kh),
\eeq

\no with $g$ the acceleration of gravity and water depth $h$. It is common to choose $c(k)=\omega_1(k)/k>0$ in \rf{whitham}, so that $\omega_1(k)$ is the root of \rf{wwdispersion} with the same sign as $k$. In what follows, we refer to this choice as the Whitham equation. The stability of periodic traveling wave solutions of the Whitham  equation has received some attention recently. Notable are \cite{ehrnstromgroveswahlen}, where the focus is on solitary waves, and \cite{hurjohnson}, where the modulational instability of small-amplitude periodic solutions is emphasized, although the main result for the high-frequency instabilities discussed in Example~\ref{ex:whitham} is included as well. Most recently, the spectral stability of periodic solutions of the Whitham equation was examined in \cite{SanfordKodamaCarterKalisch}. The goal of considering an equation like \rf{whitham} as opposed to the Korteweg-de~Vries equation or other simpler models is to capture as much of the dynamics of the full water wave problem as possible, without having to cope with the main difficulties imparted by the Euler water wave problem \cite{vandenboek} ({\em e.g.}, it is a nonlinear free boundary-value problem, the computation of its traveling-wave solutions is a nontrivial task, {\em etc}). One of the important aspects of the dynamics of a nonlinear problem is the (in)stability of its traveling wave solutions. It was shown explicitly in \cite{deconinckoliveras1} that periodic traveling wave solutions of the one-dimensional Euler water wave problem are spectrally unstable for all possible values of their parameters $h$, $g$, amplitude, and wave period. The nature of the instabilities depends on the value of these parameters. As is well known, waves in deep water are susceptible to the Benjamin-Feir or modulational instability (see \cite{zakharovostrovsky} for a review). In addition, waves in both deep and shallow water of all non-zero amplitudes are unstable with respect to high-frequency perturbations: these are perturbations whose growth rates do not have a small imaginary part, resulting in oscillatory behavior in time, independent of the spatial behavior of the unstable modes. The work in \cite{SanfordKodamaCarterKalisch} does not reveal any high-frequency instabilities for solutions of small amplitude in water that is shallow, in the context of the Whitham equation. Thus an important aspect of the Euler water wave dynamics is absent from \rf{whitham}. We will provide an analytical indication that the Whitham equation misses the presence of these instabilities, while explaining why they are missed. This explanation leads to a way to address this problem.

For suitable $N(u)$, the Whitham equation \rf{whitham} is a Hamiltonian system. In fact, for our considerations, it matters only that the linearized Whitham equation is Hamiltonian. The Lagrangian structure with the dispersion relation given by $\omega_1(k)$ as in \rf{wwdispersion} was already written down by Whitham in \cite{whitham3}, from which the Hamiltonian structure easily follows. Explicitly, for the linearized Whitham equation posed on the whole line, for any odd $\omega(k)$ we have

\beq
H=-\frac{1}{2}\int_{-\infty}^{\infty} \int_{-\infty}^\infty K(x-y) u(x)u(y)dx dy,
\eeq

\no with $J=\p_x$. Then

\beq\la{whithamham}
u_t=\p_x \frac{\delta H}{\delta u}.
\eeq

\no If instead the linearized equation is posed with periodic boundary conditions $u(x+L,t)=u(x,t)$, it follows immediately from \rf{whithamasform} that we have

\beq\la{whithamper}
u_t+\int_{-L/2}^{L/2} K(x-y)u(y)dy=0,
\eeq

\no where we have used a Fourier series instead of a Fourier transform. Further,

\beq
K(x)=\frac{1}{L}\sum_{j=-\infty}^\infty c(k_j)e^{ik_jx},
\eeq

\no and $k_j=2\pi j/L$, $j\in \mathbb{Z}$. The Hamiltonian formulation for the periodic Whitham equation \rf{whithamper} is also given by \rf{whithamham}, but with

\beq
H=-\frac{1}{2} \int_{-L/2}^{L/2}\int_{-L/2}^{L/2} K(x-y) u(x) u(y) dx dy.
\eeq

\no In fact, a formal $L\rightarrow \infty$ immediately results in the recovery of the equation posed on the whole line. Thus the Whitham equation and its periodic solutions fit in to the framework developed in this manuscript. It is one of many examples we use below. Other notable examples are the Euler water wave problem (as expected, allowing us to check our results with those of MacKay \& Saffman \cite{mackaysaffman}), the KdV equation, the Sine-Gordon equation, {\em etc}. We are particularly interested in the comparison between the results for the Euler water wave problem and those for the Whitham equation.

The results of Examples~\ref{ex:whitham} and~\ref{ex:ww} show that the Whitham equation cannot possess the high-frequency instabilities present in the water wave problem. This leads us to propose a new model equation, a so-called Boussinesq-Whitham or bidirectional Whitham equation. This equation is shown to at least satisfy the same necessary condition for the presence of high-frequency instabilities as the water wave problem, and these high-frequency instabilities originate from the same points on the imaginary axes as they do for the Euler equations. However, it is easily seen that this equation is ill posed for solutions that do not have zero average. As such it is a poor candidate to replace the Whitham equation \rf{whitham} as a shallow-water equation with the correct dispersive behavior.

\vs

{\bf Remark.} If $\omega(k)$ is not odd, then the function $u$ solving the linear Whitham equation \rf{whithamasform} is necessarily complex. The problem is Hamiltonian:

\beq\la{whithamcomplex}
u_t=-i \frac{\delta H}{\delta u^*},
\eeq

\no where $*$ denotes the complex conjugate. The Poisson structure is

\beq
J=\left(
\ba{rr}
0 & -i\\
i & 0
\ea
\right),
\eeq

\no although usually one writes only the first of the two evolution equations, omitting the equation that is the complex conjugate of \rf{whithamcomplex}. The Hamiltonian is given by

\beq
H=\int_0^{2\pi} u^* \omega(-i\p_x)u dx.
\eeq

\no A real formulation in terms of the real and imaginary parts of $u$ is possible as well (using the canonical transformation $u=(q+i p)/\sqrt{2}$), resulting in a canonical Hamiltonian structure. The linear Whitham equation \rf{whithamcomplex} is often rewritten in Fourier space using the discrete Fourier transform, due to the periodic boundary conditions. This leads to the Hamiltonian $H=\sum_{n=-\infty}^\infty \omega(n) z_n z_n^*$.

\section{Scalar Hamiltonian PDEs}\label{sec:scalar}

In this section, we investigate the stability of $2\pi$-periodic traveling wave solutions of Hamiltonian systems of the form

\beq\la{scalar}
u_t=\p_x \frac{\delta H}{\delta u},
\eeq

\no where $u(x,t)$ is a scalar real-valued function. Thus $J=\partial_x$. Since this Poisson operator is singular, all equations of this form conserve the quantity $\int_0^{2\pi} u dx$, which is the Casimir for this Poisson operator. Systems of this form include the Korteweg- de~Vries equation \cite{gardner, zakharovfaddeev} and its many generalizations, the Whitham equation \rf{whitham}, and many others. As mentioned above, our only interest is in the linearization of these equations around their trivial solution. We write the quadratic part $H^0$ of $H$ as

\beq\la{h0scalar}
H^0=-\frac{1}{2}\int_0^{2\pi} \sum_{n=0}^\infty \alpha_n u_{nx}^2 dx,
\eeq

\no where the coefficients $\alpha_n\in \mathbb{R}$. As before, indices on $u$ denote partial derivatives. Specifically, $u_{nx}$ denotes $\p_x^{n} u$.

For most examples, the number of terms in \rf{h0scalar} is finite, and all but a few of the coefficients $\alpha_n$ are nonzero. For the Whitham equation \rf{whitham}, the number of nonzero terms is infinite, but convergence is easily established. Note that \rf{h0scalar} is the most general form of a quadratic Hamiltonian depending on a single function. Indeed, a term in $H^0$ containing $u_{mx} u_{nx}$, with $m$ and $n$ positive integers, may be reduced to a squared term using integration by parts.

Using the notation \rf{h0scalar}, the linearized equation is

\beq\la{scalarlinear}
u_t=-\sum_{n=0}^\infty (-1)^n \alpha_n u_{(2n+1)x}.
\eeq

\no  We proceed with the six steps outlined in the introduction.

\begin{enumerate}

\item {\bf Quadratic Hamiltonian.} The modified Hamiltonian $H_c^0$ is given by
\beq
H_c^0=\frac{c}{2}\int_0^{2\pi}u^2dx-\frac{1}{2}\int_0^{2\pi} \sum_{n=0}^\infty \alpha_n u_{nx}^2 dx.
\eeq

\item {\bf Dispersion relation.} For equations of the form \rf{scalar}, the dispersion relation has only a single branch:

\beq\la{scalardispersion}
\omega(k)=\sum_{n=0}^\infty \alpha_n k^{2n+1}.
\eeq

The absence of even powers of $k$ in \rf{scalardispersion} is due to our imposition that \rf{scalarlinear} is a conservative equation, {\em i.e.}, there is no dissipation. All integers are allowable $k$ values, since we have equated the period to be $2\pi$. The equation \rf{scalarlinear} may be written as

\beq
u_t=-i \omega(-i\p_x)u.
\eeq

\item {\bf Bifurcation branches.} Since \rf{scalar} is scalar, only one branch can bifurcate away from the trivial solution. To find the corresponding value of $c$, we write \rf{scalarlinear} in a moving frame as

\beq\la{scalarlinearmoving}
u_t-c u_x=i \omega(i\p_x)u.
\eeq

\no This equation has its own dispersion relation given by

\beq
\Omega(k)=\omega(k)-c k,
\eeq

\no \sloppypar obtained by looking for solutions of the form $u=\exp(ikx-i \Omega t)$. Letting $u=\sum_{n=-\infty}^\infty \exp(i n x) \hat u_n$, it follows that $\p_t \hat u_n=-i \Omega(n)\hat u_n$. Thus a nonzero stationary solution may exist provided $\Omega(N)=0$, for $N\in \mathbb{N}, N\neq 0$. We have used the oddness of $\Omega(N)$ to restrict to strictly positive values of $N$. Thus the starting point of the bifurcation branch in the (speed, amplitude)-plane is $(c,0)$, where $c$ is determined by

\beq
c=\frac{\omega(N)}{N},
\eeq

\no for any integer $N>0$. Choosing $N>1$ implies that the fundamental period of the solutions is not $2\pi$, but $2\pi/N$. In practice, we choose $N=1$. A Fourier series approximation to the explicit form of the small-amplitude solutions corresponding to this bifurcation branch may be obtained using a standard Stokes expansion \cite{stokes, whitham}.

\item {\bf Stability spectrum.} In order to compute the stability spectrum associated with the zero-amplitude solution at the start of the bifurcation branch, we let $u(x,t)=U(x)\exp(\lambda t)+$c.c., where c.c. denotes the complex conjugate of the preceding term. As usual, if any $\lambda$ are found for which the real part is positive, the solution is spectrally unstable \cite{kapitulapromislow}. All bounded eigenfunctions $U(x)$ may be represented as

\beq\la{efscalar}
U(x)=\sum_{n=-\infty}^\infty a_n e^{i(n+\mu)x},
\eeq

\no where $\mu\in (-1/2,1/2]$ is the Floquet exponent. Such a representation for $U(x)$ is valid even for solutions on the bifurcation branch of nonzero amplitude \cite{deconinckkutz1}. Since \rf{scalarlinearmoving} is a problem with constant coefficients, only a single term in \rf{efscalar} is required. We obtain

\beq\la{evscalar}
\lambda_n^{(\mu)}=-i\Omega(n+\mu)=-i \omega(n+\mu)+i(n+\mu)c, ~~n\in \mathbb{Z}.
\eeq

\no As expected, all eigenvalues are imaginary and the zero-amplitude solution is neutrally stable. For a fixed value of $\mu$, \rf{evscalar} gives a point set on the imaginary axis in the complex $\lambda$ plane. As $\mu$ is varied in $(-1/2, 1/2]$, these points trace out intervals on the imaginary axis. Depending on $\omega(k)$, these intervals may cover the imaginary axis.

\item {\bf Collision condition.} The most generic scenarios for two eigenvalues given by \rf{evscalar} to collide are that (i) two of them are zero, and they collide at the origin, and (ii) two of them are equal, but nonzero. We ignore the first possibility, since the next step proves to be inconclusive for this case,  as discussed in the introduction. The second possibility requires $\lambda_n^{(\mu)}=\lambda_m^{(\mu)}$, for some $m,n\in \mathbb{Z}$, $m\neq n$, fixed $\mu\in (-1/2, 1/2]$, and $\lambda_n^{(\mu)}$, $\lambda_m^{(\mu)}\neq 0$. This may be rewritten as

    \beq\la{scalarcollision}
    \frac{\omega(n+\mu)-\omega(m+\mu)}{n-m}=\frac{\omega(N)}{N}, ~~m,n\in \mathbb{Z}, m\neq n \mbox{~and~} \mu\in (-1/2, 1/2].
    \eeq

\no This equation has an elegant graphical interpretation: the right-hand side is fixed by the choice of $N$, fixing the bifurcation branch in Step~3. It represents the slope of a line through the origin and the point $(N,\omega(N))$ in the $(k, \omega)$ plane. The left hand side is the slope of a line in the same plane passing through the points $(n+\mu, \omega(n+\mu))$ and $(m+\mu, \omega(m+\mu))$, see Fig.~\ref{fig:scalarslope}.

Even though the graph of the dispersion relation admits parallel secant lines, this is not sufficient for a solution of \rf{scalarcollision}, as it is required that their abscissas are an integer apart. Nevertheless, the graphical interpretation can provide good intuition for solving the collision condition, which typically has to be done numerically.

%
\begin{figure}[tb]
\def\svgwidth{5.8in}
\centerline{\hspace*{0in}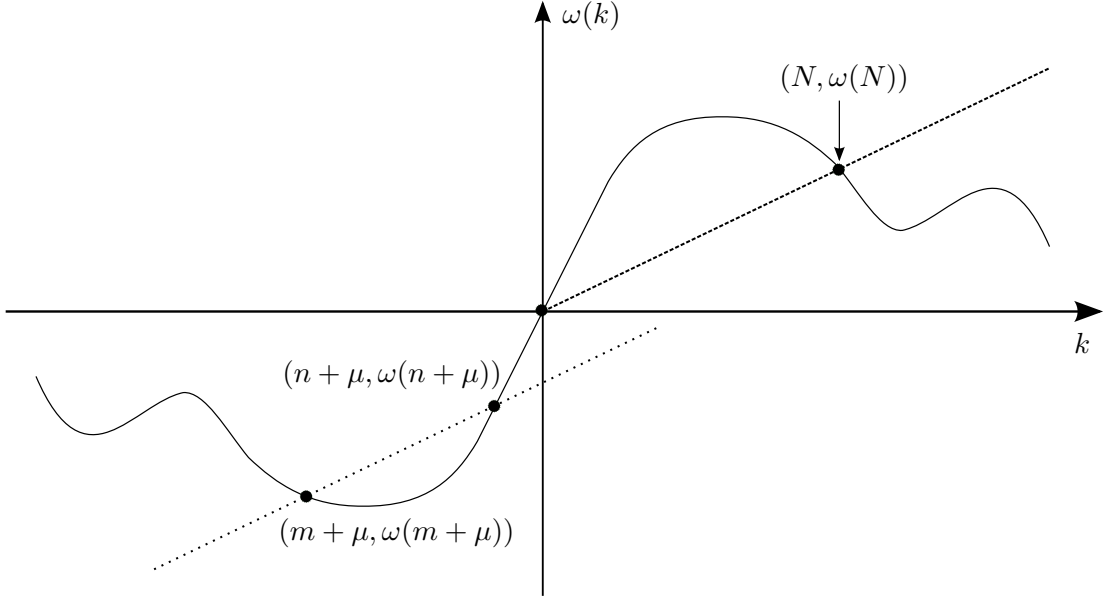}
\caption{\la{fig:scalarslope} The graphical interpretation of the collision condition \rf{scalarcollision}. The solid curve is the graph of the dispersion relation $\omega(k)$. The slope of the dashed line in the first quadrant is the right-hand side in \rf{scalarcollision}. The slope of the parallel dotted line is its left-hand side.
    }
\end{figure}
%

\item {\bf Krein signature.} The Krein signature of an eigenvalue is the sign of the Hamiltonian $H^0_c$ evaluated on the eigenspace associated with the eigenvalue. We are considering two simple eigenvalues colliding, thus the eigenspace for each eigenvalue consists of multiples of the eigenfunction only. To allow for eigenfunctions of the form $a_n \exp(i(n+\mu)x+\lambda_n^{(\mu)}t)+$c.c, which are not $2\pi$-periodic (unless $\mu=0$), it is necessary to replace the integral in \rf{h0scalar} with a whole-line average. More details on this process are found, for instance, in \cite{deconinckoliveras1}. A simple calculation shows that the contribution to $H^0_c$ from the $(n,\mu)$ mode is proportional to the ratio of $\Omega(n+\mu)/(n+\mu)$:

    \beq\la{scalarkrein}
    H^0_c|_{(n,\mu)}\sim -|a_n|^2 \frac{\Omega(n+\mu)}{n+\mu}.
    \eeq

\no Other terms are present in the Hamiltonian density, but they have zero average. The sign of this expression is the Krein signature of the eigenvalue $\lambda_n^{(\mu)}$. Thus a necessary condition for $\lambda_n^{(\mu)}$ and $\lambda_m^{(\mu)}$ to leave the imaginary axis for solutions of non-zero amplitude is that the signs of \rf{scalarkrein} with $(n,\mu)$ and $(m, \mu)$ are different, contingent on $\mu$, $m$ and $n$ satisfying \rf{scalarcollision}. Explicitly, this condition is

\beq\la{scalarkreincondition1}
\sgn\left[\frac{\omega(n+\mu)}{n+\mu}-c\right]\neq \sgn\left[\frac{\omega(m+\mu)}{m+\mu}-c\right].
\eeq

\no Alternatively, the product of the left-hand side and the right-hand side should be negative. Using \rf{scalarcollision}, \rf{scalarkreincondition1} becomes

\beq\la{scalarkreincondition2}
(n+\mu)(m+\mu)<0,
\eeq

\no or, provided $mn\neq 0$, and using that $\mu \in (-1/2, 1/2]$,

\beq\la{scalarkreincondition3}
n m<0.
\eeq

\end{enumerate}

{\bf Remarks.}

\begin{itemize}

\item It is clear from \rf{scalarkrein} why our methods do not lead to any conclusions about collisions of eigenvalues at the origin. If $\lambda_n^{(\mu)}=0$, then $\Omega(n+\mu)=0$, and the contribution to the Hamiltonian of such a mode vanishes. As a consequence, the associated Krein signature is zero.

\item When the theory of \cite{kollarmiller} is restricted to the case of solutions of zero-amplitude, so as to recover the constant coefficient stability problem, the graphical stability criterion given there coincides with the one presented here.

\end{itemize}

\vs

We conclude our general considerations of this section with the following summary.

\vs

Assume that the linearization of the scalar Hamiltonian system \rf{scalar} is dispersive ({\em i.e.,} its dispersion relation $\omega(k)$ is real valued for $k\in \mathbb{R}$). Let $N$ be a strictly positive integer.
Consider $2\pi/N$-periodic traveling wave solutions of this system of sufficiently small-amplitude and with velocity sufficiently close to $\omega(N)/N$. In order for these solutions to be spectrally unstable with respect to high-frequency instabilities as a consequence of two-eigenvalue collisions, it is necessary that there exist $n$, $m~\in \mathbb{Z}$, $n\neq m$, $\mu \in (-1/2, 1/2]$ for which

\beq
\frac{\omega(n+\mu)}{n+\mu}\neq \frac{\omega(N)}{N}, ~~\frac{\omega(m+\mu)}{m+\mu}\neq \frac{\omega(N)}{N},
\eeq

\no such that

\beq
\frac{\omega(n+\mu)-\omega(m+\mu)}{n-m}=\frac{\omega(N)}{N},
\eeq

\no and

\beq
(m+\mu)(n+\mu) <0.
\eeq

Next we proceed with some examples.

\subsection{The (generalized) Korteweg - de Vries equation}

We consider the generalized KdV (gKdV) equation

\beq\la{gkdv}
u_t+\sigma u^n u_x +u_{xxx}=0,
\eeq

\no where we restrict $n$ to integers 1 or greater. Here $\sigma$ is a constant coefficient, chosen as convenient. Important special cases discussed below are the KdV equation ($n=1$) and the modified KdV (mKdV) equation ($n=2$). Many of the details below extend easily to more general nonlinearities, with the main requirement being that the linearized equation is $u_t+u_{xxx}=0$. The stability of periodic solutions of the gKdV equation has received some attention recently \cite{bronskijohnson1, bronskijohnsonkapitula, todd, johnsonzumbrunbronski}. For the integrable cases $n=1$ and $n=2$, more detailed analysis is possible, see \cite{bottmandeconinck, deconinckkapitula, deconincknivalamkdv}. We do not claim to add anything new to these discussions, but we wish to use this example to illustrate how the six-step process outlined in this section leads to easy conclusions before moving on to more complicated settings.

\begin{enumerate}

\item The modified Hamiltonian is given by

\beq
H^0_c=\frac{1}{2}\int_0^{2\pi} (u_x^2+c u^2)dx.
\eeq

\item The dispersion relation is

\beq
\omega=-k^3.
\eeq

\item Bifurcation branches in the ($c$, amplitude)-plane start at $(c,0)$, with

\beq
c=\frac{\omega(k)}{k}=-k^2.
\eeq

\no Since we desire $2\pi$ periodic solutions, we choose $k=1$. Any choice $k=N$, where $N$ is a non-zero integer is allowed. Choosing $k=1$, bifurcation branches start at $(-1,0)$.

For the integrable cases $n=1$ (KdV) and $n=2$ (mKdV), these bifurcation branches may be calculated in closed form. For the KdV equation in a frame traveling with speed $c$, the $2\pi$-periodic solutions are given by (with $\sigma=1$)

\beq\la{cnkdv}
u=\frac{12\kappa^2 K^2(\kappa)}{\pi^2}\cn^2\left(\frac{K(\kappa)x}{\pi},\kappa\right),
\eeq

\no where $\cn$ denotes the Jacobian elliptic cosine function and $K(\kappa)$ is the complete elliptic integral of the first kind \cite{dlmf, dlmfbook}. Further,

\beq\la{ckdv}
c(\kappa)=\frac{4K^2(\kappa)}{\pi^2}(2\kappa^2-1),
\eeq

\no resulting in an explicit bifurcation curve $(c(\kappa), 12\kappa^2 K^2(\kappa)/\pi^2$), parameterized by the elliptic modulus $\kappa\in [0,1)$.   This bifurcation curve is shown in Fig.~\ref{fig:kdvbranch}a.

%
\begin{figure}[tb]
\begin{tabular}{cc}
\multicolumn{2}{c}{\def\svgwidth{3.5in}
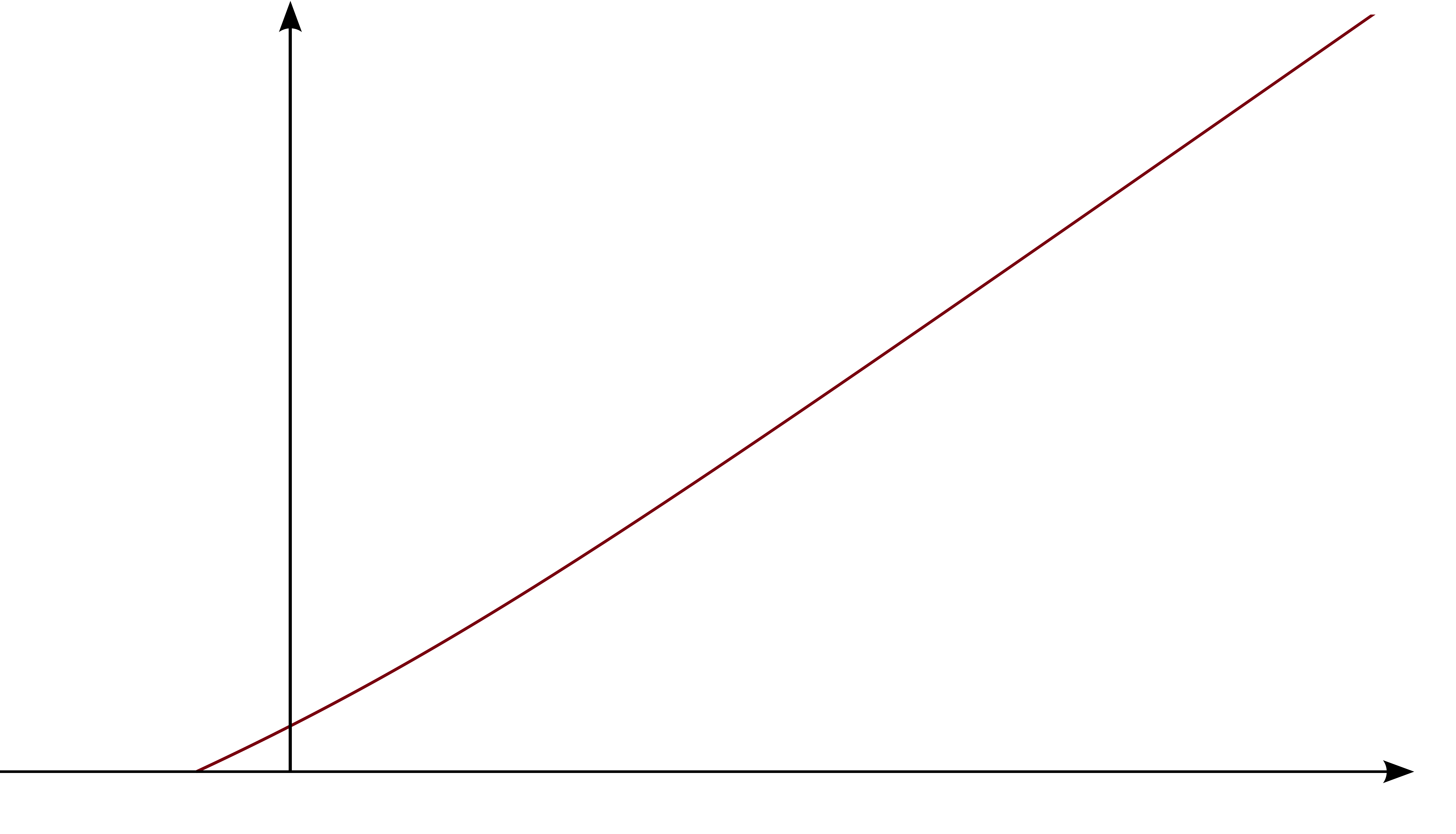}\\
\multicolumn{2}{c}{(a)}\\
~&~\\
\def\svgwidth{2.8in}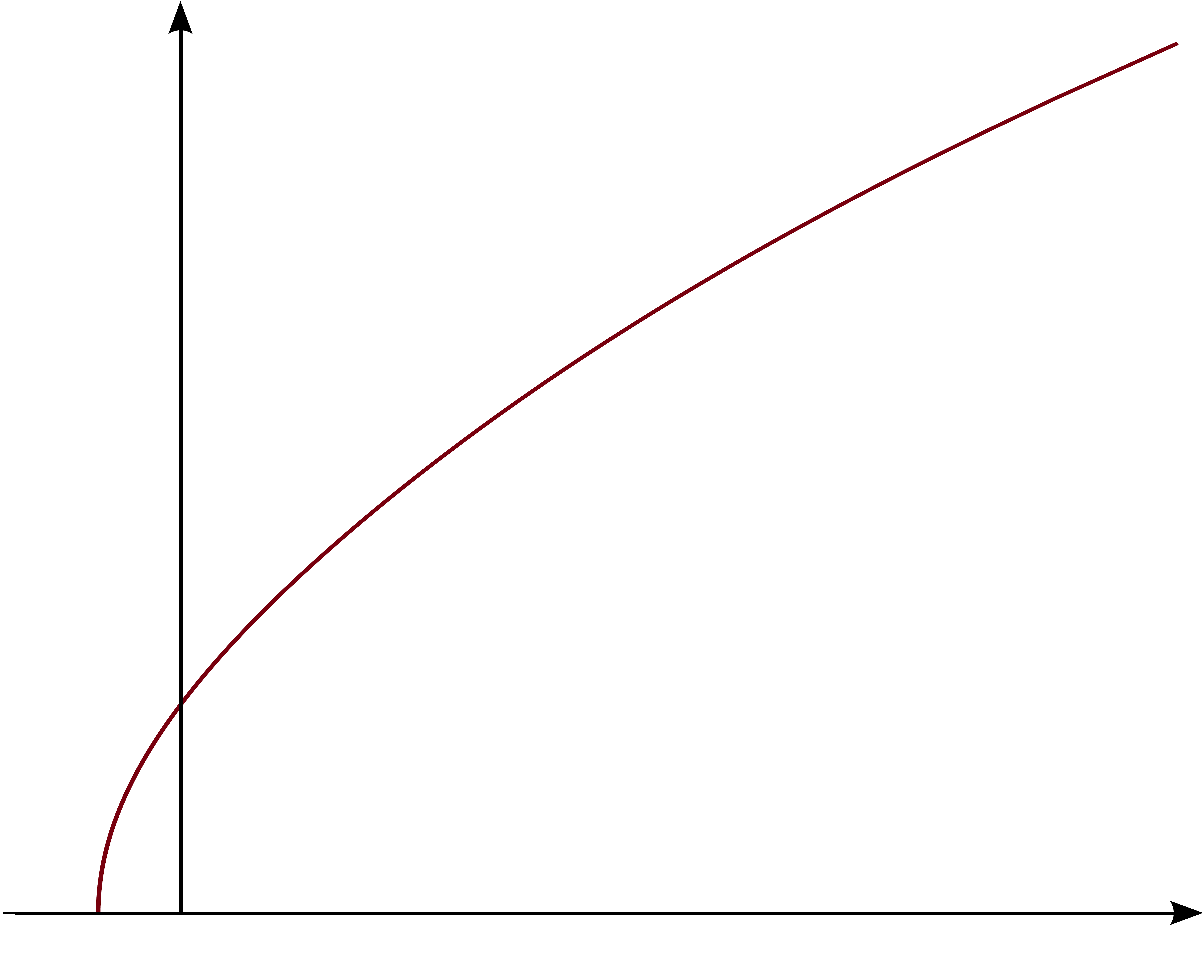 &
\def\svgwidth{2.8in}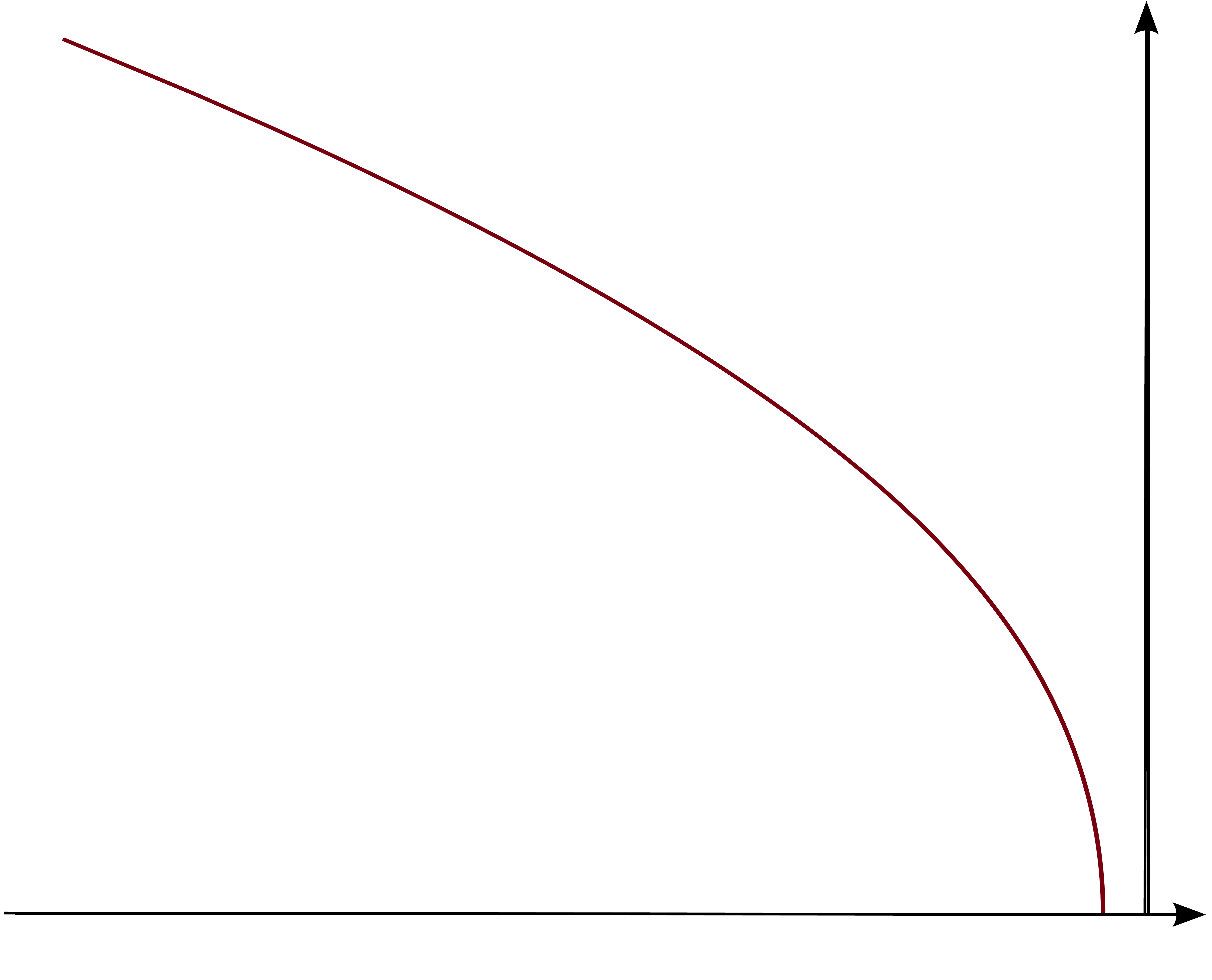\\
(b) & (c)
\end{tabular}
\caption{\la{fig:kdvbranch} The amplitude {\em vs.} $c$ bifurcation plots for the traveling-wave solutions of the generalized KdV equation \rf{gkdv}. (a) The KdV equation, $n=1$, for the cnoidal wave solutions \rf{cnkdv}. (b) The mKdV equation, $n=2$, for the cnoidal wave solutions \rf{mkdvcn}. Lastly, (c) shows the bifurcation plot for the snoidal wave solutions \rf{mkdvsn} of mKdV, $n=2$.
Note that all bifurcation branches start at $(-1,0)$, as stated above.
    }
\end{figure}
%

For the mKdV equation ($n=2$), different families of traveling-wave solutions exist \cite{deconincknivalamkdv}. We consider two of the simplest. For $\sigma=3$ (focusing mKdV), a family of $2\pi$-periodic solutions is given by

\beq\la{mkdvcn}
u=\frac{2\sqrt{2}\kappa K(\kappa)}{\pi}\cn\left(\frac{2 K(\kappa)x}{\pi},\kappa\right),
\eeq

\no with $c(\kappa)$ given by \rf{ckdv}, resulting in an explicit bifurcation curve $(c(\kappa), 2\sqrt{2}\kappa K(\kappa)/\pi)$, parameterized by the elliptic modulus $\kappa\in [0,1)$.   This bifurcation curve is shown in Fig.~\ref{fig:kdvbranch}b. It should be noted that a solution branch exists where the solution is expressed in terms of the Jacobian $\dn$ function \cite{dlmf, dlmfbook}: $u=\sqrt(2)K(\kappa)\dn(K(\kappa)x/\pi,\kappa)/\pi$, but this solution does not have a small-amplitude limit and our methods do not apply directly to it. Rather, the solutions limit to the constant solution $u=1/\sqrt{2}$ as $\kappa\ra 0$. A simple transformation $v=u-1/\sqrt{2}$ transforms the problem to one where our methods apply.

For $\sigma=-3$ (defocusing mKdV), a period $2\pi$ solution family is

\beq\la{mkdvsn}
u=\frac{2\sqrt{2}\kappa K(\kappa)}{\pi}\sn\left(\frac{2 K(\kappa)x}{\pi},\kappa\right),
\eeq

\no with $c(\kappa)=-4(1+\kappa^2)K^2(\kappa)/\pi^2$. Here $\sn$ is the Jacobian elliptic sine function \cite{dlmf, dlmfbook}, resulting in an explicit bifurcation curve $(c(\kappa), 2\sqrt{2}\kappa K(\kappa)/\pi)$, parameterized by the elliptic modulus $\kappa\in [0,1)$.   This bifurcation curve is shown in Fig.~\ref{fig:kdvbranch}c.

\item The stability spectrum is given by \rf{evscalar}, with $\omega(k)=-k^3$ and $c=-1$, resulting in

\beq
\lambda^{(n)}_\mu=i(n+\mu)(1+(n+\mu)^2).
\eeq

These eigenvalues cover the imaginary axis, as $n$ and $\mu$ are varied. The imaginary part of this expression is displayed in Fig.~\ref{fig:threecover}a. For the sake of comparison with Fig.~2 in \cite{bottmandeconinck} we let $\mu\in [-1/4, 1/4)$, which implies that $n$ is any half integer. The results of Fig.~2 in \cite{bottmandeconinck} are for elliptic modulus $\kappa=0.8$, implying a solution of moderate amplitude. The comparison of these two figures serves to add credence to the relevance of the results obtained using the zero-amplitude solutions at the start of the bifurcation branch.

%
\begin{figure}[tb]
\begin{tabular}{cc}
\def\svgwidth{2.8in}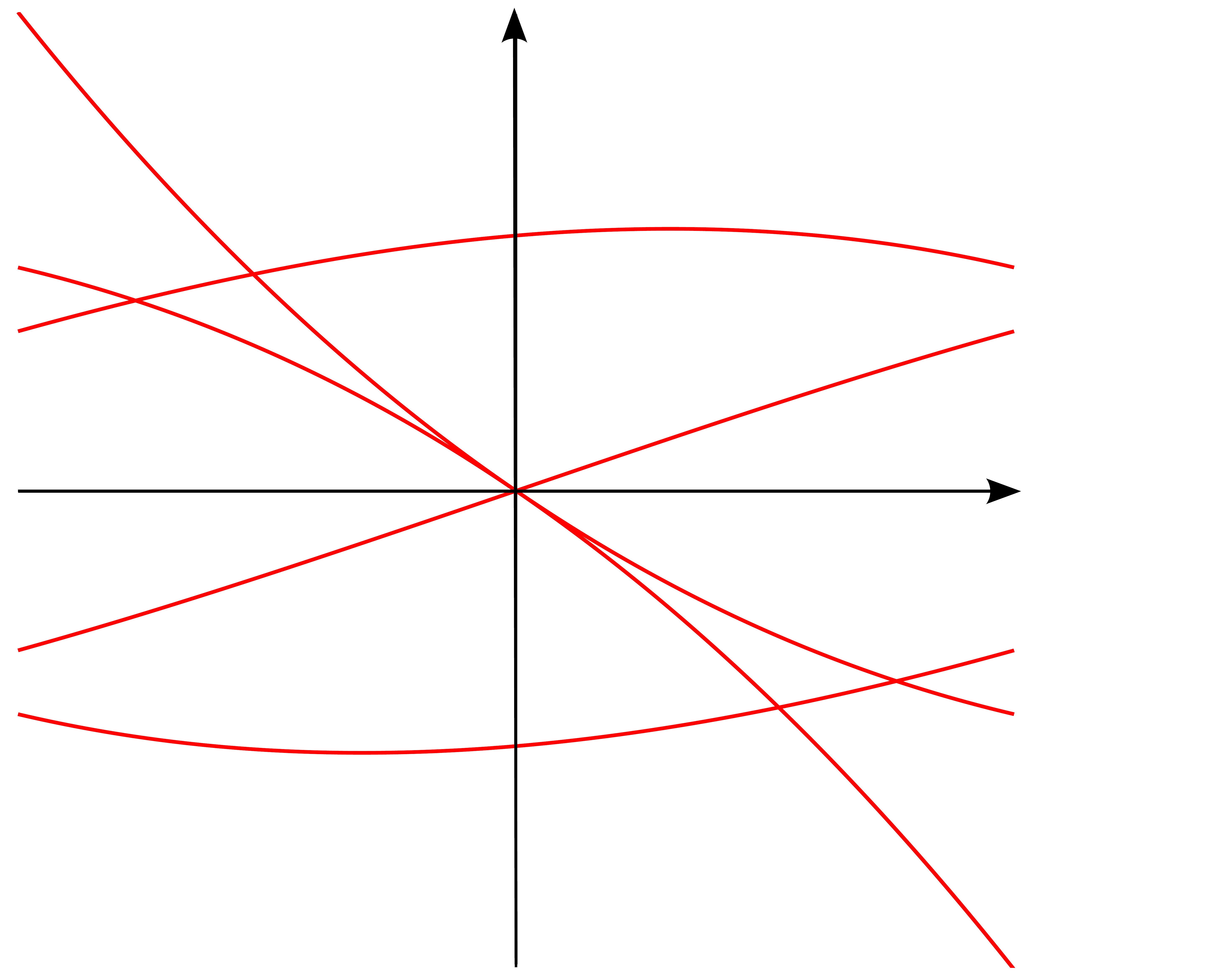 &
\def\svgwidth{2.8in}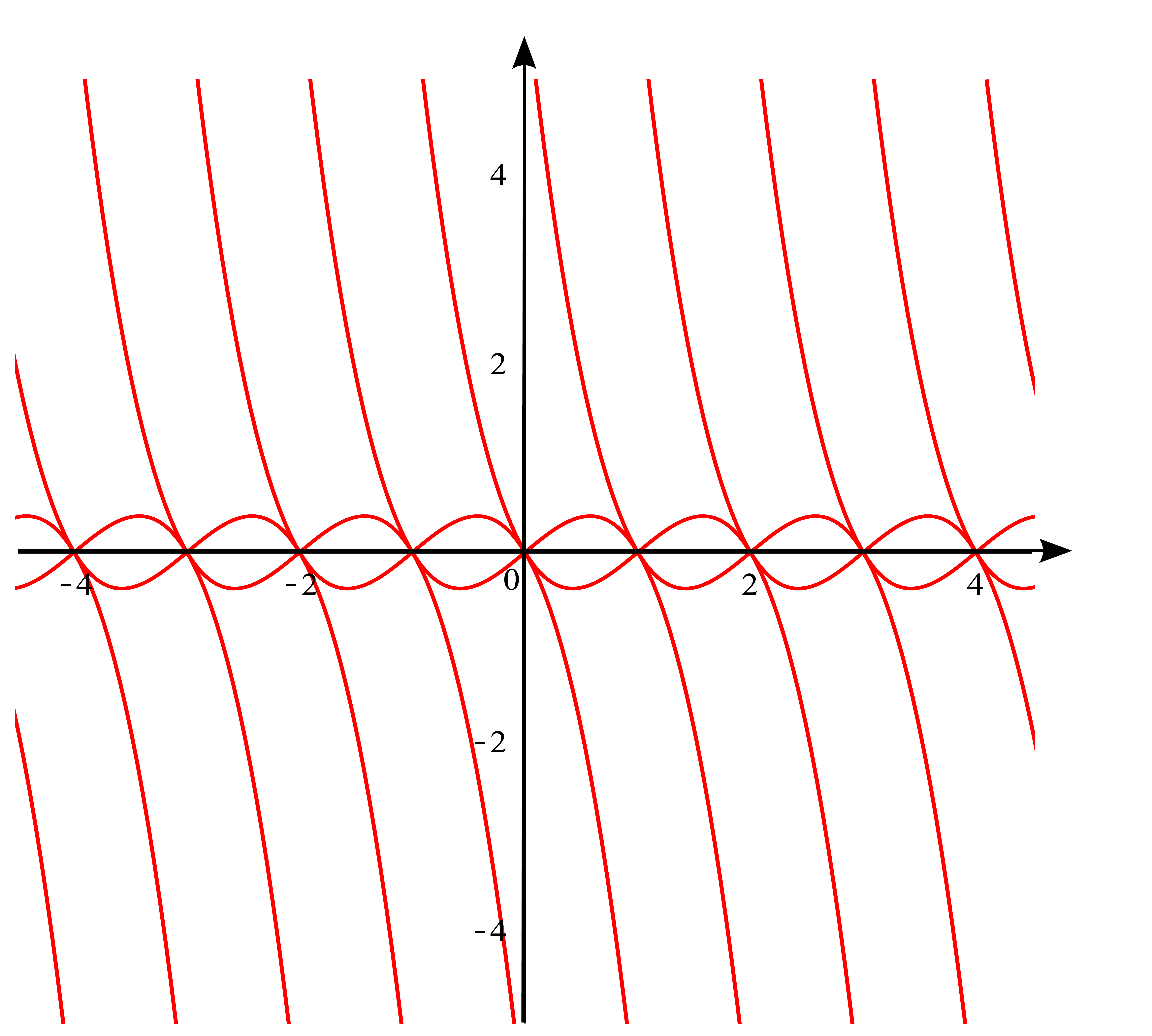\\
(a) & (b)
\end{tabular}
\caption{\la{fig:threecover} (a) The imaginary part of $\lambda_n^{(\mu)}\in (-0.7, 0.7)$ as a function of $\mu\in[-1/4, 1/4)$. (b) The curves $\Omega(k+n)$, for various (integer) values of $n$, illustrating that collisions occur at the origin only.
    }
\end{figure}
%

\item With $n+\mu=k$ and $m+\mu=k+l$, for some $l\in \mathbb{Z}$, the collision condition \rf{scalarcollision} is written as

\beq
l^2+3kl+3k^2-1=0,
\eeq

\no where the trivial solution $l=0$ has been discarded. This is the equation of an ellipsoid in the $(k,l)$ plane. It intersects lines of nonzero integer $l$ in six integer points: $\pm (1, -2)$, $\pm(0, 1)$, $\pm(1, -1)$. Since for all of these, $\Omega(k)=0$, any collisions happen only at the origin $\lambda_n^{(\mu)}=0$. This is also illustrated in Fig.~\ref{fig:threecover}b.

\item The final step of our process is preempted by the results of the previous step. No Krein signature of colliding eigenvalues can be computed, since no eigenvalues collide.

\end{enumerate}

Since eigenvalues do not collide away from the origin they cannot leave the imaginary axis through such collisions and no high-frequency instabilities occur for small amplitude solutions of the gKdV equation. This result applies to the KdV and mKdV equations as special cases. The absence of high-frequency instabilities for small amplitude solutions is consistent with the results in, for instance, \cite{bottmandeconinck, deconinckkapitula, todd}.

\subsection{The Whitham equation}\la{ex:whitham}

As our second scalar example, we consider the Whitham equation \rf{whitham}. For this example, no analytical results exist, but the work of Sanford {\em et al.} \cite{SanfordKodamaCarterKalisch} allows for a comparison with numerical results. Sanford {\em et al.} do not report the presence of high-frequency instabilities for solutions of any period. Their absence has been verified by us using the same methods, see Fig.~\ref{fig:wstuff}. Hur \& Johnson \cite{hurjohnson} consider periodic solutions, focusing on the modulational instability. However, they do include a Krein signature calculation of the eigenvalues of the zero-amplitude solutions, reaching the same conclusions obtained below. In what follows the nonlinear term $N(u)$ does not contribute, as in the previous examples.

%
\begin{figure}[tb]
\begin{tabular}{cc}
\def\svgwidth{3.2in}\hspace*{-0.2in}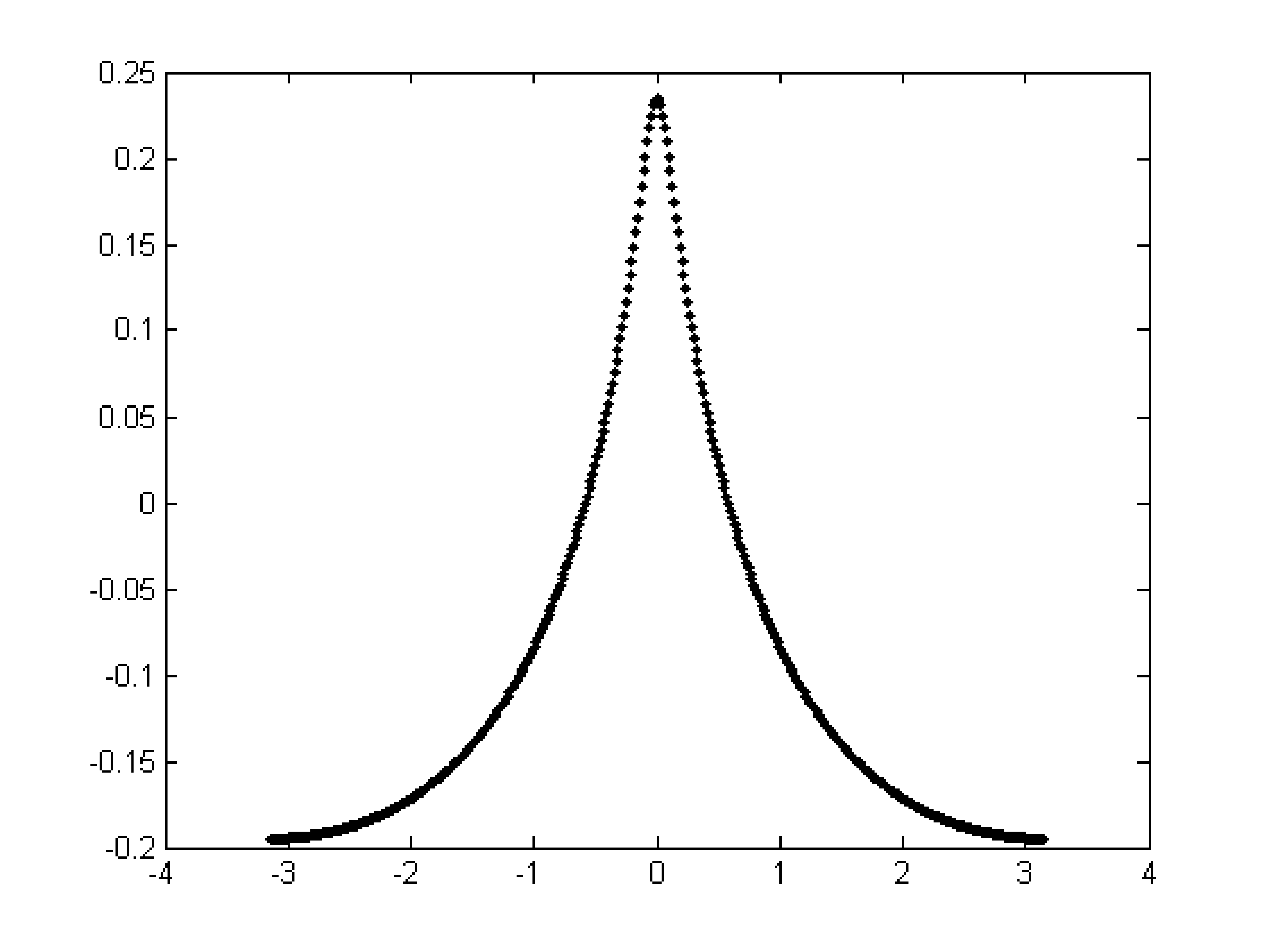 &
\def\svgwidth{3.6in}\hspace*{-0.4in}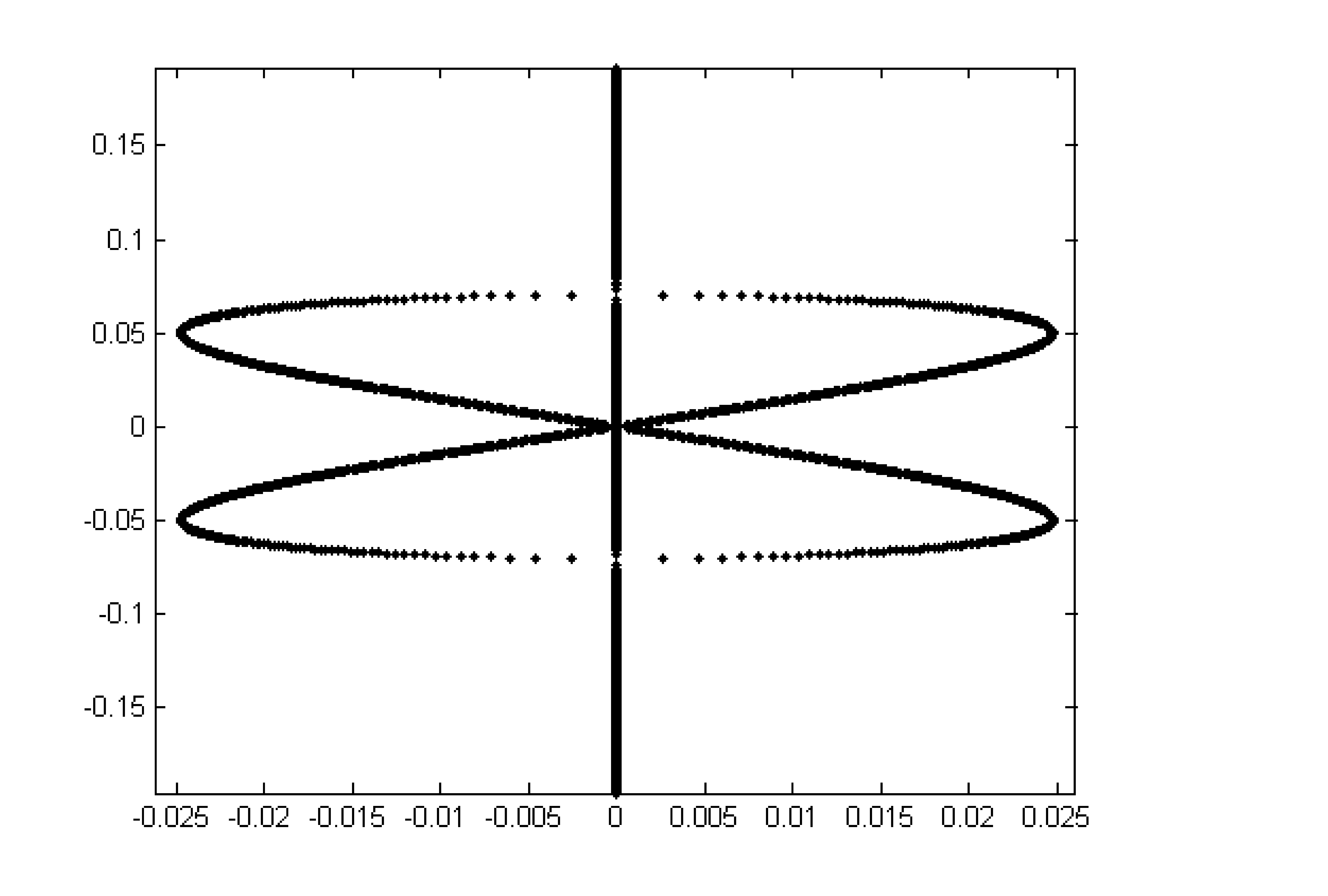\\
(a) & \hspace*{-0.4in}(b)
\end{tabular}
\caption{\la{fig:wstuff} (a) The profile of a $2\pi$-periodic small-amplitude traveling wave solution of the Whitham equation \rf{whitham} with $c\approx 0.7697166847$, computed using a cosine collocation method with 128 Fourier modes, see \cite{SanfordKodamaCarterKalisch}. (b) The stability spectrum of this solution, computed using the Fourier-Floquet-Hill method \cite{deconinckkutz1} with $128$ modes and 2000 different values of the Floquet parameter $\mu$. The presence of a modulational instability is clear, but no high-frequency instabilities are observed, in agreement with the theory presented. Note that the hallmark bubbles of instability were looked for far outside of the region displayed here.
    }
\end{figure}
%

\begin{enumerate}

\item The modified Hamiltonian is $$H^0_V=\frac{V}{2}\int_0^{2\pi} u^2 dx-\frac{1}{2}\int_{-\pi}^\pi \int_{-\pi}^\pi K(x-y)u(x) u(y) dx dy.$$

\no We use $V$ to denote the speed of the traveling wave, to avoid confusion with the phase speed $c(k)$ in the kernel of the Whitham equation.

\item The dispersion relation is given by $\omega(k)=\sgn(k)\sqrt{gk\tanh(kh)}$.

\item The bifurcation branch starts at $(V,0)=(\sqrt{g \tanh(h)},0)$, where we have chosen $N=1$ so that the minimal period of the solutions is $2\pi$.

\item The elements of the stability spectrum are given by

\beq\la{evwhitham}
\lambda_n^{(\mu)}=i(n+\mu)\sqrt{g \tanh(h)}-i \sgn(n+\mu)\sqrt{g(n+\mu)\tanh (h(n+\mu))}.
\eeq

\item The dispersion relation for the Whitham equation is plotted in Fig.~\ref{fig:wcollision}(a), together with the line through the origin with slope $\omega(1)/1$. Since the dispersion relation is concave down (up) in the first (third) quadrant, the condition \rf{scalarcollision} is not satisfied. Thus collisions of eigenvalues away from the
    origin do not occur. This is also illustrated in Fig.~\ref{fig:wcollision}(b), where the imaginary part of $\lambda_n^{(\mu)}$ is plotted for various integers $n$.

%
\begin{figure}[tb]
\begin{tabular}{cc}
\def\svgwidth{2.8in}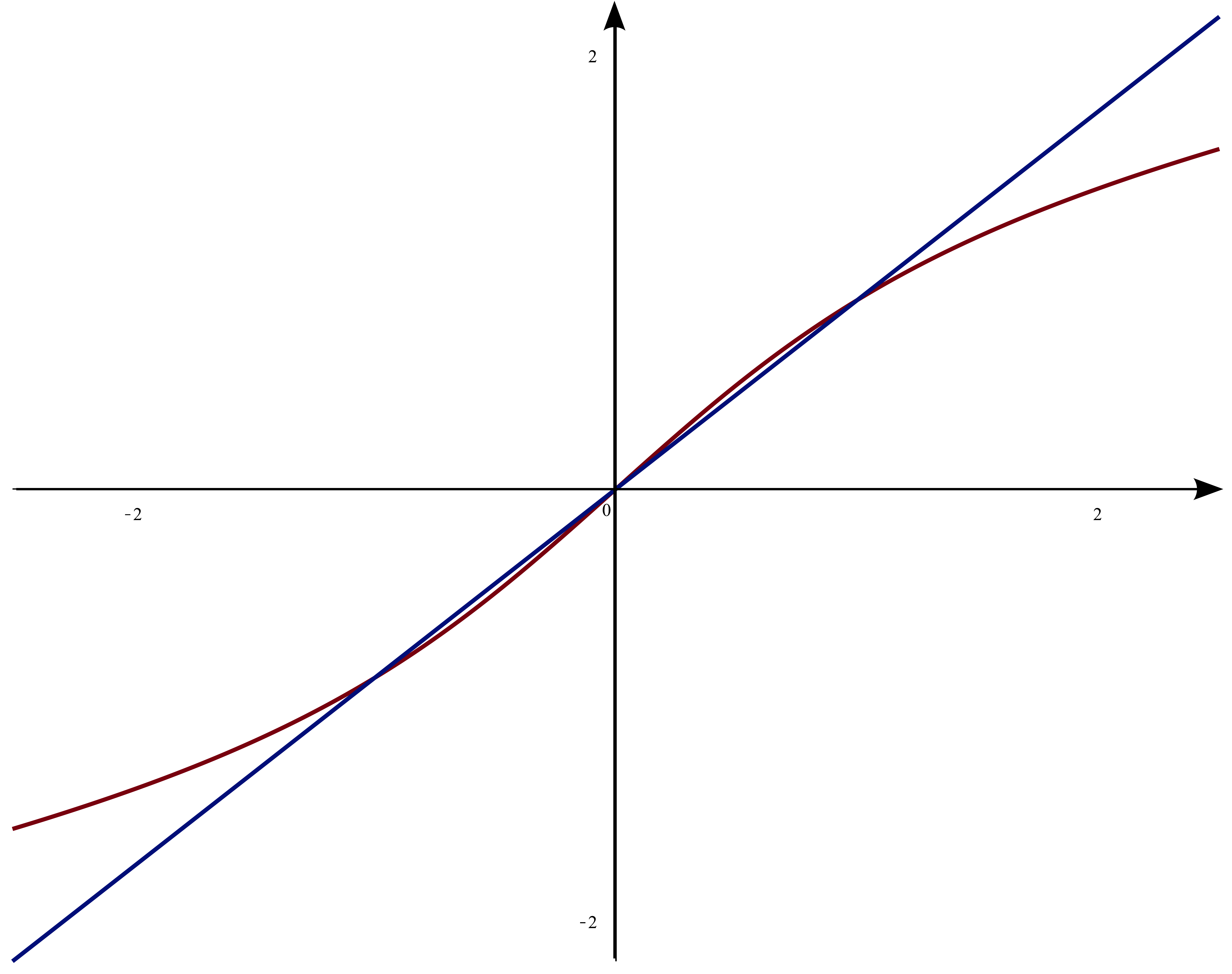 &
\def\svgwidth{2.8in}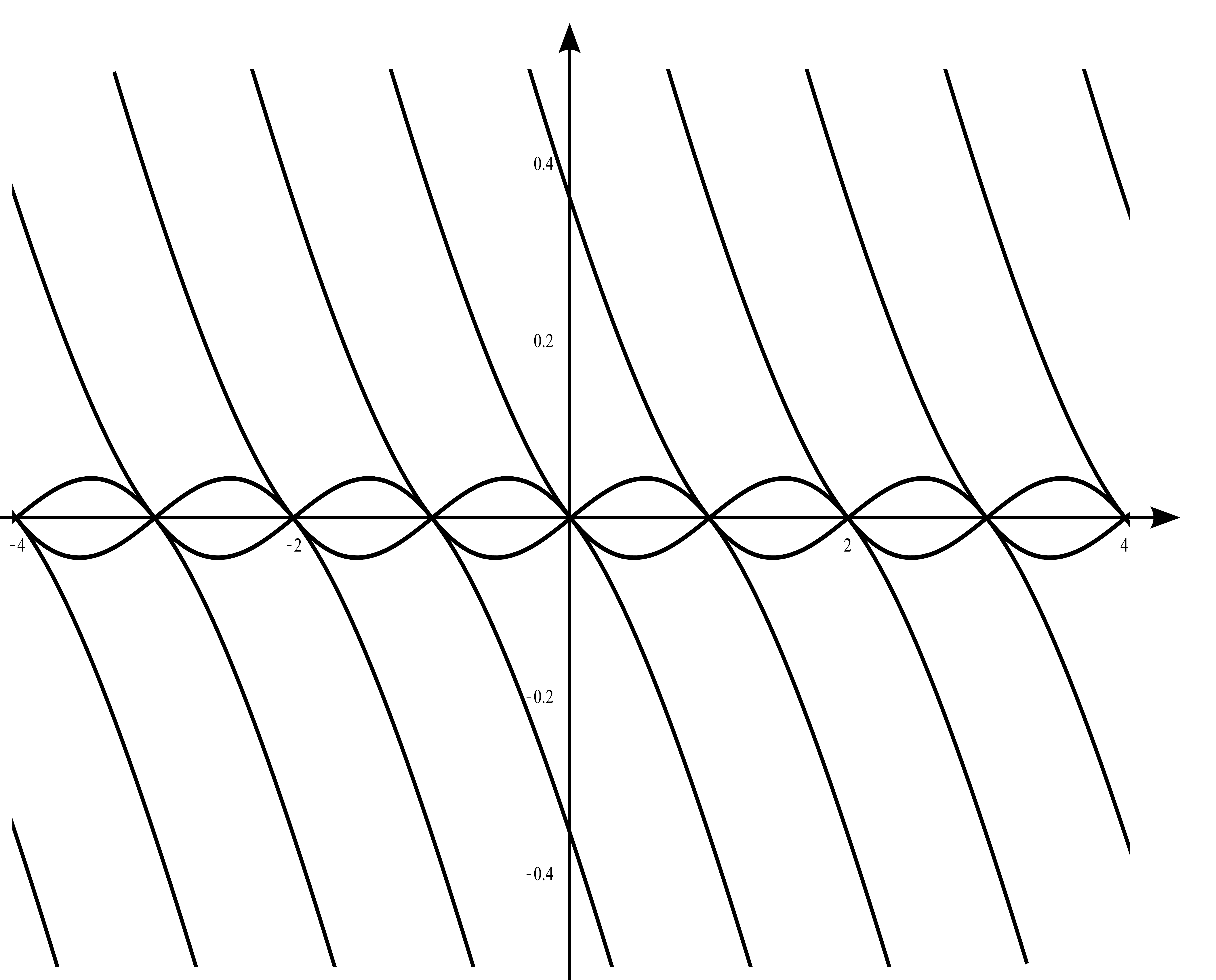\\
(a) & (b)
\end{tabular}
\caption{\la{fig:wcollision} (a) The dispersion relation for the Whitham equation (curve), together with the line through the origin of slope $\omega(1)/1$, representing the right-hand side of \rf{scalarcollision}.  (b) The curves $\Omega(k+n)$, for various (integer) values of $n$, illustrating that collisions occur at the origin only.
    }
\end{figure}
%

\item No Krein signature calculation is relevant since eigenvalues do not collide away from the origin.

\end{enumerate}

We conclude that periodic solutions of sufficiently small amplitude of the Whitham equation are not susceptible to high-frequency instabilities. This is consistent with the results presented in \cite{SanfordKodamaCarterKalisch}, see also Fig.~\ref{fig:wstuff}. Thus, the Whitham equation is unable to replicate the instabilities found in the shallow depth water wave problem for solutions of small amplitude, despite having a dispersion relation that is identical to one branch of the water wave dispersion relation. We return to this in the next section.

\section{Two-component Hamiltonian PDEs with canonical Poisson structure}\la{sec:vector}

We generalize the ideas of the previous section to the setting of two-component Hamiltonian PDEs with canonical Poisson structure. In other words, the evolution PDE can be written as

\begin{align}\la{can}
\pp{}{t}
\left(
\ba{c}
q\\p
\ea
\right)&=J \nabla H ~~\iff~~\left\{
\ba{rcl}
\ds q_t&=&\ds \dd{H}{p}\\
\ds p_t&=&\ds -\dd{H}{q}
\ea
\right.,
\end{align}

\no where the Poisson operator $J$ is given by \rf{canonicalJ}. This Poisson operator is nonsingular, thus there are no Casimirs. Examples of systems of this form are the Nonlinear Schr\"odinger equation in real coordinates \cite{bottmandeconincknivala}, the Sine-Gordon equation \cite{faddeev, jonesmiller1, jonesmiller2}, and the water wave problem \cite{zakharov1}. As before, our interest is in the linearization of this system around the zero-amplitude solution. The quadratic Hamiltonian corresponding to this linearization can be written as

\beq\la{canham}
H^0=\int_0^{2\pi} \left(
\frac{1}{2}\sum_{n=0}^\infty c_n q_{nx}^2+\frac{1}{2}\sum_{n=0}^\infty b_n p_{nx}^2+p\sum_{n=0}^\infty a_n q_{nx}
\right)dx,
\eeq

\no with $a_n$, $b_n$, $c_n\in \mathbb{R}$. Typically the number of terms in the sums above is finite, but an example like the water wave problem requires the possibility of an infinite number of nonzero contributing terms in the Hamiltonian. As for the Whitham equation, convergence of the resulting series is not problematic. The form \rf{canham} is the most general form of a quadratic Hamiltonian depending on two functions $q(x,t)$ and $p(x,t)$. Indeed, any quadratic term of a form not included above is reduced to a term that is included by straightforward integration by parts. The linearization of \rf{can} is given by

\alpheqn

\begin{align}\la{canlin}
q_t&=\sum_{n=0}^\infty a_n q_{nx}+\sum_{n=0}^\infty (-1)^n b_n p_{2nx},\\
p_t&=-\sum_{n=0}^\infty (-1)^n c_n q_{2nx}-\sum_{n=0}^\infty (-1)^n a_n p_{nx}.
\end{align}

\resetalpheqn

We proceed with the six step program outlined in the introduction.

\begin{enumerate}

\item {\bf Quadratic Hamiltonian.} The modified Hamiltonian $H_c^0$ is given by

\beq
H_c^0=\int_0^{2\pi} \left(c p q_x+
\frac{1}{2}\sum_{n=0}^\infty c_n q_{nx}^2+\frac{1}{2}\sum_{n=0}^\infty b_n p_{nx}^2+p\sum_{n=0}^\infty a_n q_{nx}
\right)dx.
\eeq

\no This expression serves as a repository for the coefficients which are needed in what follows.

\item {\bf Dispersion relation.} We look for solutions to \rf{canlin} of the form $q=\hat q \exp(ikx-i\omega t)$, $p=\hat p \exp(ikx-i\omega t)$. Requiring the existence of non-trivial ({\em i.e.}, non-zero) solutions, we find that $\omega(k)$ is determined by

\beq\la{candis}
\det \left(
\ba{cc}
\ds i\omega+\sum_{n=0}^\infty a_n (ik)^n &
\ds \sum_{n=0}^\infty b_n k^{2n}\\
\ds -\sum_{n=0}^\infty c_n k^{2n} &
\ds i\omega -\sum_{n=0}^\infty a_n (-1)^n (ik)^n
\ea
\right)=0.
\eeq

\no This is a quadratic equation for $\omega(k)$, resulting in two branches of the dispersion relation, $\omega_1(k)$ and $\omega_2(k)$. Assuming that (\ref{canlin}-b) is dispersive, $\omega_1(k)$ and $\omega_2(k)$ are real-valued for $k\in \mathbb{R}$. This is not easily translated in a condition on the coefficients $a_n$, $b_n$, $c_n$ and $d_n$, since their reality is not assumed.

\item {\bf Bifurcation branches.} Traveling wave solutions are stationary solutions of

\alpheqn

\begin{align}\la{linmoving}
q_t&=c q_x+\sum_{n=0}^\infty a_n q_{nx}+\sum_{n=0}^\infty (-1)^n b_n p_{2nx}=\dd{H_c^0}{p},\\
p_t&=c p_x-\sum_{n=0}^\infty (-1)^n c_n q_{2nx}-\sum_{n=0}^\infty (-1)^n a_n p_{nx}=-\dd{H_c^0}{q}.
\end{align}

\resetalpheqn

\no This system has the dispersion relations $\Omega_{1,2}(k)=\omega_{1,2}(k)-c k$. In Fourier space the stationary equations become

\alpheqn

\begin{align}
0&=ik c \hat q+\sum_{n=0}^\infty a_n (ik^n) \hat q+\sum_{n=0}^\infty (-1)^n b_n (ik)^{2n} \hat p,\\
0&=ik c \hat p-\sum_{n=0}^\infty (-1)^n c_n (ik)^{2n} \hat q-\sum_{n=0}^\infty (-1)^n a_n (ik)^n \hat p.
\end{align}

\resetalpheqn

\no Thus $c$ is obtained from the condition that these equations have a nontrivial solution $(\hat q, \hat p)$. This condition requires that the $2\times 2$ determinant of the system above is zero. A simple comparison with \rf{candis} gives that there are two bifurcation points given by $(\omega_1(N)/N,0)$ and $(\omega_2(N)/N,0)$. Any positive integer value of $N$ is allowed, but we usually choose $N=1$ so that the fundamental period is $2\pi$. In what follows, we examine the small-amplitude solutions starting from the branch $(c,0)=(\omega_1(N)/N,0)$, without loss of generality.

In many systems the two bifurcation branches are reflections of each other about the vertical axis. The corresponding solution profiles are identical to each other, moving to the right on one branch, moving to the left on the other. Examples are given below. Non-symmetric bifurcation branches cannot be excluded, however, without imposing extra assumptions on the coefficients of (\ref{canlin}-b).

\item {\bf Stability spectrum.} To find the stability spectrum, we let $q=Q(x) \exp(\lambda t)$, $p=P(x) \exp(\lambda t)$. Next, using Floquet's Theorem,

\beq\la{canfloquet}
Q=e^{i\mu x}\sum_{j=-\infty}^\infty Q_j e^{i j x}, ~~P=e^{i\mu x}\sum_{j=-\infty}^\infty P_j e^{ijx},
\eeq

\no with $\mu\in (-1/2, 1/2]$. Since (\ref{linmoving}-b) has constant coefficients, it suffices to consider monochromatic waves, {\em i.e.}, only one term of the sums in \rf{canfloquet} is retained. It follows that $\lambda$ satisfies \rf{candis} with $i\omega$ replaced by $-\lambda+i(n+\mu)c$. Thus

\beq\la{canlambda}
\lambda_{n,l}^{(\mu)}=i(n+\mu)c-i\omega_l(n+\mu)=-i \Omega_l(n+\mu), ~~~~l=1,2.
\eeq

\no As expected, the zero solution is neutrally stable since $\omega_{1,2}(k)$ are real, assuming dispersive equations. The stability spectrum consists of two one-parameter point sets, one for $l=1$, the other for $l=2$.

\item {\bf Collision conditions.} Ignoring collisions at the origin, we require $\lambda_{n_1, l_1}^{(\mu)}=\lambda_{n_2, l_2}^{(\mu)}\neq 0$ for some $n_1, n_2 \in \mathbb{Z}$, $\mu\in (-1/2,1/2]$, $l_1, l_2\in \{1,2\}$. This gives

\beq\la{cancollision}
\frac{\omega_{l_1}(n_1+\mu)-\omega_{l_2}(n_2+\mu)}{n_1-n_2}=\frac{\omega_1(N)}{N}.
\eeq

\no The right-hand side depends on $\omega_1$ since we have chosen the first branch of the dispersion relation in Step~3. As before, this collision condition may be interpreted as a parallel secant condition, but with the additional freedom of being able to use points from both branches of the dispersion relation. This is illustrated in Fig.~\ref{fig:cancollision}.

%
\begin{figure}[tb]
\def\svgwidth{4.8in}
\centerline{\hspace*{0in}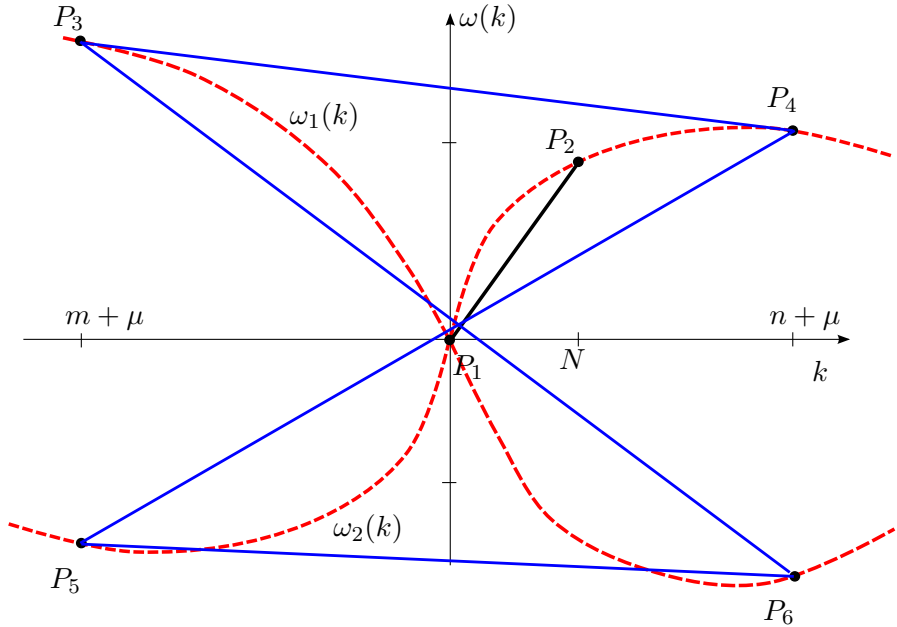}
\caption{\la{fig:cancollision} The graphical interpretation of the collision condition \rf{cancollision}. The dashed curves are the graphs of the dispersion relations $\omega_1(k)$ and $\omega_2(k)$. The slope of the segment $P_1P_2$ is the right-hand side in \rf{cancollision}. The collision condition \rf{cancollision} seeks points whose abscissas are an integer apart, so that at least one of the segments $P_3P_4$, $P_3P_6$, $P_5P_4$ or $P_5P_6$ is parallel to the segment $P_1P_2$.
    }
\end{figure}
%

\item {\bf Krein signature.} In the setting of a system of Hamiltonian PDEs as opposed to a scalar PDE, we use a different but equivalent characterization of the Krein signature~\cite{mackay}. The Krein signature is the contribution to the Hamiltonian of the mode involved with the collision. Since our Hamiltonians are quadratic, this implies that the Krein signature of the eigenvalue $\lambda$ with eigenvector $v$ is given by

    \beq\la{kreinhessian}
    \mbox{signature}(\lambda, v)=\sgn(v^\dagger {\cal L}_c v),
    \eeq

\no where ${\cal L}_c$ is the Hessian of the Hamiltonian $H_c^0$, and $v^\dagger$ denotes the complex conjugate of the transposed vector. Since the Hessian ${\cal L}_c$ is a symmetric linear operator, the argument of the $\sgn$ in \rf{kreinhessian} is real. Recall that the linearization of the system \rf{hamc} can be written as

\beq
\p_t\left(
\ba{c}
q\\p
\ea
\right)=
J {\cal L}_c
\left(
\ba{c}
q\\p
\ea
\right),
\eeq

\no which makes it easy to read off ${\cal L}_c$. For the case of (\ref{linmoving}-b),

\beq\la{hessian}
{\cal L}_c=\left(
\ba{cc}
\ds \sum_{n=0}^\infty c_n (-1)^n \p_x^{2n} &
\ds -c \p_x+\sum_{n=0}^\infty a_n (-1)^n \p_x^n \\
\ds c \p_x+\sum_{n=0}^\infty a_n \p_x^n &
\ds \sum_{n=0}^\infty b_n (-1)^n \p_x^{2n}
\ea
\right).
\eeq

\no Next, the eigenvectors $v$ are given by

\beq
\left(
\ba{c}
q \\ p
\ea
\right)=e^{\lambda t+i\mu x+i n x}\left(
\ba{c}
Q_n \\ P_n
\ea
\right),
\eeq

\no where $(Q_n, P_n)^T$ satisfies

\beq\la{evec}
\lambda
\left(
\ba{c}
\ds Q_n \\ \ds P_n
\ea
\right)=J
\hat {\cal L}_c\left(\ba{c}
\ds Q_n \\ \ds P_n
\ea
\right).
\eeq

\no Here $\hat {\cal L}_c$ is the symbol of ${\cal L}_c$, {\em i.e.}, the $2\times 2$ matrix obtained by replacing $\p_x\ra i(n+\mu)$ in \rf{hessian}:

\beq\la{symbol}
\hat {\cal L}_c=
\left(
\ba{cc}
\ds \sum_{n=0}^\infty c_n (n+\mu)^{2n} &
\ds -i c (n+\mu)+\sum_{n=0}^\infty a_n (-1)^n (in+i\mu)^n \\
\ds i c (n+\mu) +\sum_{n=0}^\infty a_n (in+i\mu)^n &
\ds \sum_{n=0}^\infty b_n (n+\mu)^{2n}
\ea
\right).
\eeq

\no Using \rf{evec}, \rf{kreinhessian} is rewritten as

\beq
\mbox{signature}\left(\lambda_{n,l}^{(\mu)}, v_{n,l}^{(\mu)}\right)=\sgn\left(\lambda_{n,l}^{(\mu)}\det \left(
\ba{cc}
Q_n & P_n\\
Q_n^* & P_n^*
\ea
\right)\right).
\eeq

\no The determinant is imaginary, since interchanging the rows gives the complex conjugate result. Since $\lambda_{n,l}^{(\mu)}$ is imaginary, the result is real and the signature is well defined. Again, it is clear that no conclusions can be drawn if $\lambda_{n,l}^{(\mu)}=0$. Since we wish to examine whether signatures are equal or opposite, we consider the product of the signatures corresponding to $\lambda_{n_1,l_1}^{(\mu)}$ and $\lambda_{n_2,l_2}^{(\mu)}$. Using \rf{symbol} we find that signatures are opposite, provided that

\begin{align}\nonumber
&\sum_{j_1=0}^\infty c_{j_1} (n_1+\mu)^{2j_1}\sum_{j_2=0}^\infty c_{j_2} (n_2+\mu)^{2j_2}
\left(\omega_{l_1}(n_1+\mu)+\sum_{j_3=0}^\infty a_{2j_3+1}(-1)^{j_3}(n_1+\mu)^{2j_3+1}\right)\times\\\la{cankrein1}
&~~~~~~~~~~~~~~~~\left(\omega_{l_2}(n_2+\mu)+\sum_{j_4=0}^\infty a_{2j_4+1}(-1)^{j_4}(n_2+\mu)^{2j_4+1}\right)<0.
\end{align}

\no The above condition is obtained by expressing the eigenvectors in \rf{evec} in terms of the entries of the first row of \rf{symbol}. An equivalent condition is obtained using the second row:

\begin{align}\nonumber
&\sum_{j_1=0}^\infty b_{j_1} (n_1+\mu)^{2j_1}\sum_{j_2=0}^\infty b_{j_2} (n_2+\mu)^{2j_2}
\left(\omega_{l_1}(n_1+\mu)-\sum_{j_3=0}^\infty a_{2j_3+1}(-1)^{j_3}(n_1+\mu)^{2j_3+1}\right)\times\\\la{cankrein2}
&~~~~~~~~~~~~~~~~\left(\omega_{l_2}(n_2+\mu)-\sum_{j_4=0}^\infty a_{2j_4+1}(-1)^{j_4}(n_2+\mu)^{2j_4+1}\right)<0.
\end{align}

\no Depending on the system at hand, the condition \rf{cankrein1} or \rf{cankrein2} may be more convenient to use.

\end{enumerate}

{\bf Remark.} An important class of systems is those for which $\omega_1(k)=-\omega_2(k)$. We refer to such systems as even systems. It follows immediately from \rf{candis} that for even systems $a_{2j+1}=0$, $j=1,2,\ldots$. The Krein conditions \rf{cankrein1} and \rf{cankrein2} simplify significantly, becoming

\beq\la{sym1}
\omega_{l_1}(n_1+\mu)\omega_{l_2}(n_2+\mu)\sum_{j_1=0}^\infty c_{j_1} (n_1+\mu)^{2j_1}\sum_{j_2=0}^\infty c_{j_2} (n_2+\mu)^{2j_2}<0,
\eeq

\no or

\beq\la{sym2}
\omega_{l_1}(n_1+\mu)\omega_{l_2}(n_2+\mu)\sum_{j_1=0}^\infty b_{j_1} (n_1+\mu)^{2j_1}\sum_{j_2=0}^\infty b_{j_2} (n_2+\mu)^{2j_2}<0.
\eeq

We summarize our results.

\vs

Assume that the linearization of the Hamiltonian system \rf{can} is dispersive ({\em i.e.,} its dispersion relations $\omega_1(k)$ and $\omega_2(k)$ are real valued for $k\in \mathbb{R}$). Let $N$ be a strictly positive integer.
Consider $2\pi/N$-periodic traveling wave solutions of this system of sufficiently small-amplitude and with velocity sufficiently close to $\omega_1(N)/N$. In order for these solutions to be spectrally unstable with respect to high-frequency instabilities as a consequence of two-eigenvalue collisions, it is necessary that there exist $l_1$, $l_2~\in \{1,2\}$, $n_1$, $n_2~\in \mathbb{Z}$, $n_1\neq n_2$, $\mu \in (-1/2, 1/2]$ for which

\beq
\frac{\omega_{l_1}(n_1+\mu)}{n_1+\mu}\neq \frac{\omega_1(N)}{N}, ~~\frac{\omega_{l_2}(n_2+\mu)}{n_2+\mu}\neq \frac{\omega(N)}{N},
\eeq

\no such that

\beq
\frac{\omega_{l_1}(n_1+\mu)-\omega_{l_2}(n_2+\mu)}{n_1-n_2}=\frac{\omega_1(N)}{N},
\eeq

\no and \rf{cankrein1}, or equivalently, \rf{cankrein2} holds.

\vs

We proceed with examples.

\subsection{The Sine-Gordon equation}

As a first example, we consider the Sine-Gordon (SG) equation \cite{scott}:

\beq
u_{tt}-u_{xx}+\sin u=0.
\eeq

\no The stability of the periodic traveling wave solutions of this equation has been studied recently by Jones {\em et al.} \cite{jonesmiller1, jonesmiller2}. Different classes of periodic traveling wave solutions exist, but only two of those can be considered as small-amplitude perturbations of a constant background state. We consider the so-called superluminal ($c^2>1$) librational waves. The subluminal ($c^2<1$) librational waves require the use of the transformation $v=u-\pi$ so that their small amplitude limit approaches the zero solution. We do not consider them here. The limits of the rotational waves are either soliton solutions or have increasingly larger amplitude. As such the rotational waves do not fit in the framework of this paper. An overview of the properties of these solutions as well as illuminating phase-plane plots are found in \cite{jonesmiller1}. In contrast to \cite{jonesmiller1, jonesmiller2}, we fix the period of our solutions, as elsewhere in this paper. This makes a comparison of the results more complicated.

\begin{enumerate}

\item {\bf Quadratic Hamiltonian.} With $q=u$, $p=u_t$,

\beq
H_c^0=\int_0^{2\pi} \left(
c p q_x+\frac{1}{2}p^2+\frac{1}{2}q^2+\frac{1}{2}q_x^2
\right)dx.
\eeq

\no Thus $b_0=1$, $c_0=1$, $c_1=1$ are the only non-zero coefficients.

\item {\bf Dispersion relation.} Using \rf{candis},

\beq
\omega_{1,2}=\pm \sqrt{1+k^2}.
\eeq

\no These expressions are real valued for $k\in \mathbb{R}$, thus the SG equation is dispersive when linearized around the superluminal librational waves. Both branches of the dispersion relation are displayed in Fig.~\ref{fig:sgcollision}a.

\item {\bf Bifurcation branches.} With $N=1$, we obtain $c=\omega_1(1)/1=\sqrt{2}$.

\item {\bf Stability spectrum.} The stability spectrum is given by \rf{canlambda}:

\beq
\lambda_{n,l}^{(\mu)}=-i\Omega_l(n+\mu)=i(n+\mu)\sqrt{2}\mp i\sqrt{1+(n+\mu)^2},
\eeq

\no with $l=1$ ($l=2$) corresponding to the $-$ ($+$) sign. Here $n\in \mathbb{Z}$, $\mu\in [-1/2, 1/2)$.

\item {\bf Collision condition.} The collision condition \rf{cancollision} becomes

\beq
\frac{\omega_1(n_1+\mu)-\omega_2(n_2+\mu)}{n_1-n_2}=\sqrt{2}.
\eeq

\no We have chosen $\omega_{l_1}=\omega_1$ and $\omega_{l_2}=\omega_2$, since it is clear that the collision condition can only be satisfied if points from both dispersion relation branches are used. This is illustrated in Fig.~\ref{fig:sgcollision}a. In fact, many collisions occur, as is illustrated in Fig.~\ref{fig:sgcollision}b. One explicit solution is given by

\beq
n_1=3,~n_2=0,~\mu=\frac{\sqrt{10}-3}{2}\approx 0.081138830.
\eeq

%
\begin{figure}[tb]
\begin{tabular}{cc}
\def\svgwidth{2.8in}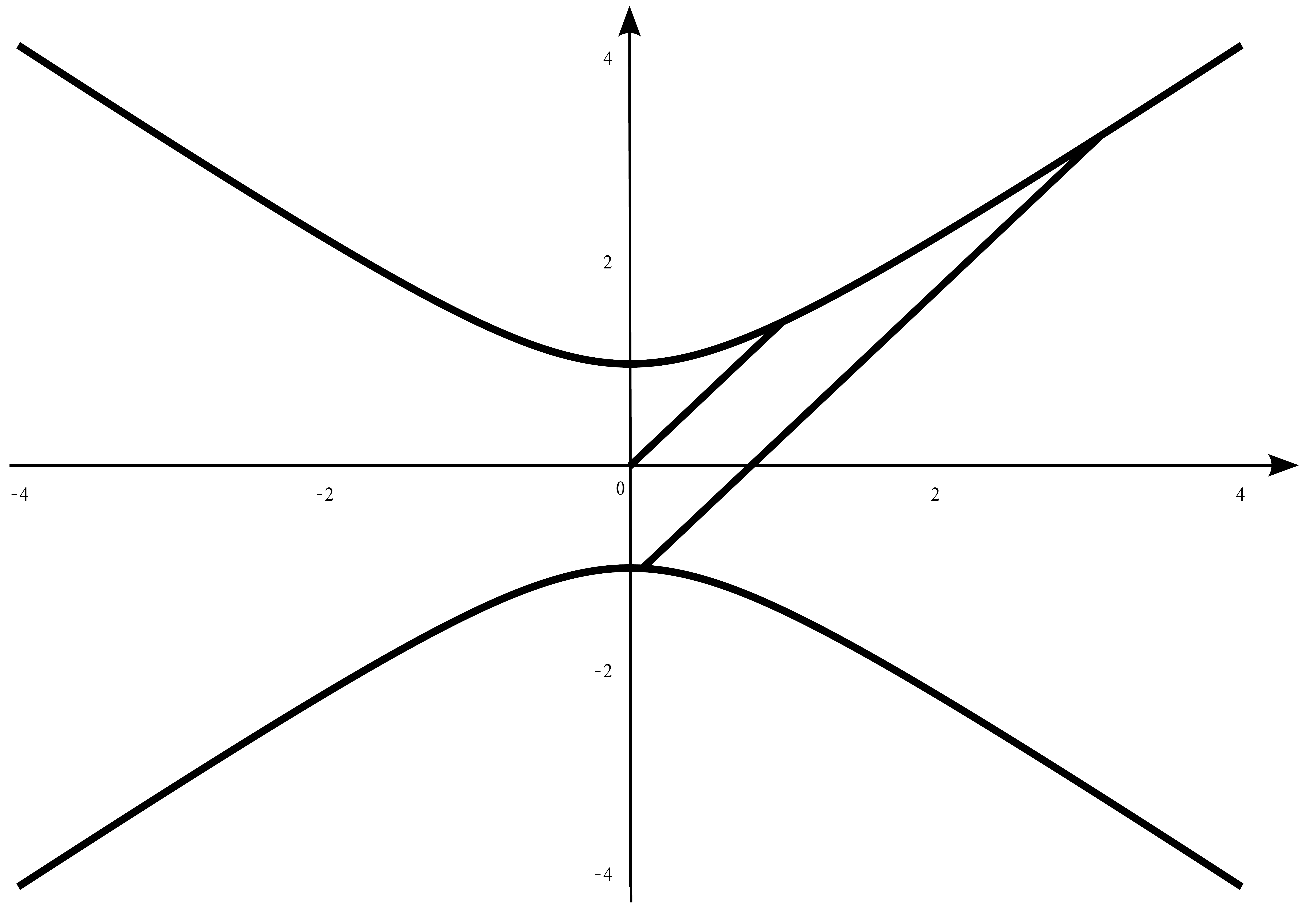 &
\def\svgwidth{2.8in}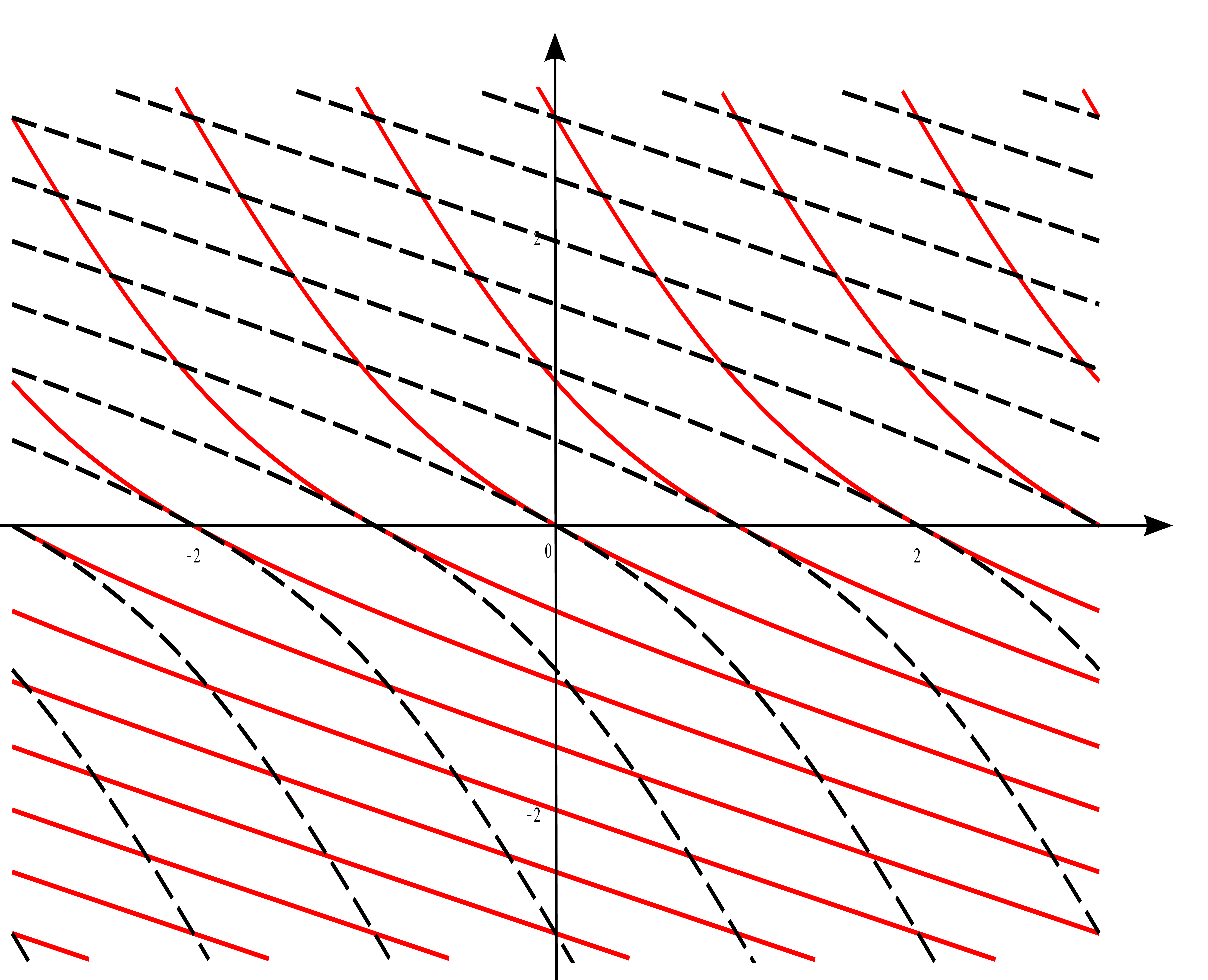\\
(a) & (b)
\end{tabular}
\caption{\la{fig:sgcollision} (a) The two branches of the dispersion relation for the Sine-Gordon equation. The line segment $P_1 P_2$ has slope $\omega(1)/1$, representing the right-hand side of \rf{cancollision}. The slope of the parallel line segment $P_3 P_4$ represents the left-hand side of \rf{cancollision}.  (b) The two families of curves $\Omega_1(k+n)$ (red, solid) and $\Omega_2(k+n)$ (black, dashed), for various (integer) values of $n$, illustrating that many collisions occur away from the origin.
    }
\end{figure}
%

\item {\bf Krein signature.} Since $\omega_2(k)=-\omega_1(k)$, we may use the conditions \rf{sym1} or \rf{sym2}. Since only one $b_j\neq 0$, \rf{sym2} is (slightly) simpler to use. We get that

\beq
\omega_1(n_1+\mu)\omega_2(n_2+\mu)<0
\eeq

\no is a necessary condition for the presence of high-frequency instabilities of small-amplitude superluminal librational solutions of the SG equation. The condition is trivially satisfied as it was remarked in the previous step that points from both dispersion relation branches have to be used to have collisions.

\end{enumerate}

%
\begin{figure}[tb]
\begin{tabular}{cc}
\def\svgwidth{3.2in}\hspace*{-0.2in}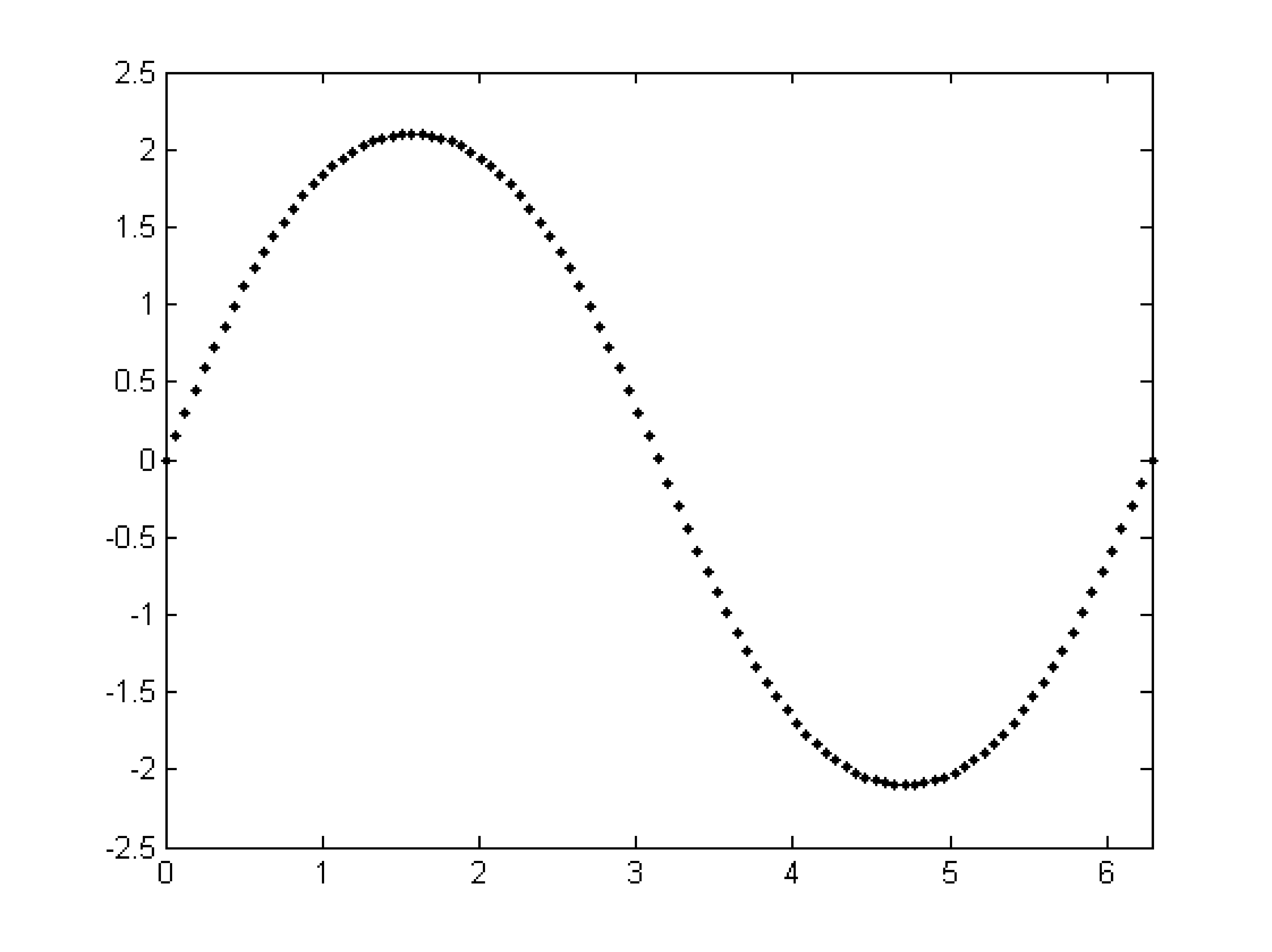 &
\def\svgwidth{3.7in}\hspace*{-0.4in}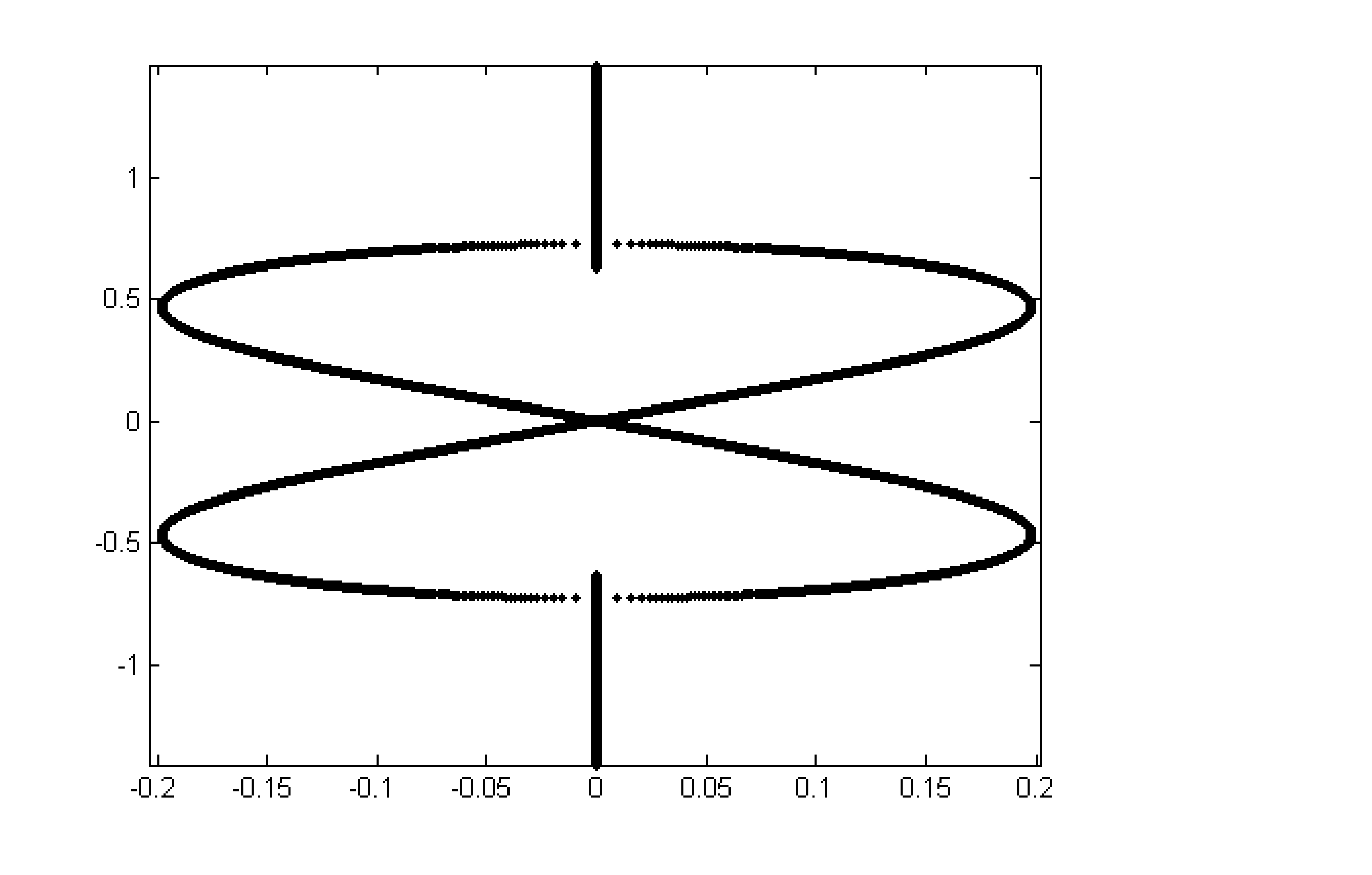\\
(a) & \hspace*{-0.4in}(b)
\end{tabular}
\caption{\la{fig:sgstuff} (a) A small-amplitude $2\pi$-periodic superluminal solution of the SG equation ($c\approx 1.236084655663$). (b) A blow-up of the numerically computed stability spectrum in a neighborhood of the origin, illustrating the presence of a modulational instability, but the absence of high-frequency instabilities.
    }
\end{figure}
%

It follows that for the superluminal solutions of the SG equations the necessary condition for the occurrence of high-frequency instabilities is satisfied. Nevertheless, as the results of \cite{jonesmiller1, jonesmiller2} show, such instabilities do not occur. This is illustrated in Fig.~\ref{fig:sgstuff}. The left panel illustrates an exact $2\pi$-periodic superluminal solution of the SG equation with $c\approx 1.236084655663$, obtained using elliptic functions. Computing the stability spectrum (right panel) of the solution using the Fourier-Floquet-Hill method \cite{deconinckkutz1} with $51$ Fourier modes and $1000$ Floquet exponents shows that no high-frequency instabilities are present to within the accuracy of the numerical method. This is consistent with the results of \cite{jonesmiller1, jonesmiller2} where only the presence of a modulational instability is observed. Thus the example of this section illustrates that the necessary condition is not always sufficient.

\subsection{The water wave problem}\la{ex:ww}

As a final example, we consider the water wave problem: the problem of determining the dynamics of the surface of an incompressible, irrotational fluid under the influence of gravity. For this example, the effects of surface tension are ignored and we consider only two-dimensional fluids, {\em i.e.}, the surface is one dimensional. The Euler equations governing the dynamics are

\alpheqn

\begin{align}\la{eulera}
\phi_{xx}+\phi_{zz}&=0, && (x,z)\in D,\\
\phi_z&=0, && z=-h,\\
\eta_t+\eta_x \phi_x&=\phi_z,&& z=\eta(x,t),\\
\phi_t+\frac{1}{2}(\phi_x^2+\phi_z)+g\eta&=0, && z=\eta(x,t),
\end{align}

\resetalpheqn

\no where $x$ and $z$ are the horizontal and vertical coordinate, respectively, see Fig.~\ref{fig:ww}; $z=\eta(x,t)$ is the free top boundary and $\phi(x,z,t)$ is the velocity potential. Further, $g$ is the acceleration due to gravity and $h$ is the average depth of the fluid.

%
\begin{figure}[tb]
\centerline{\def\svgwidth{6in}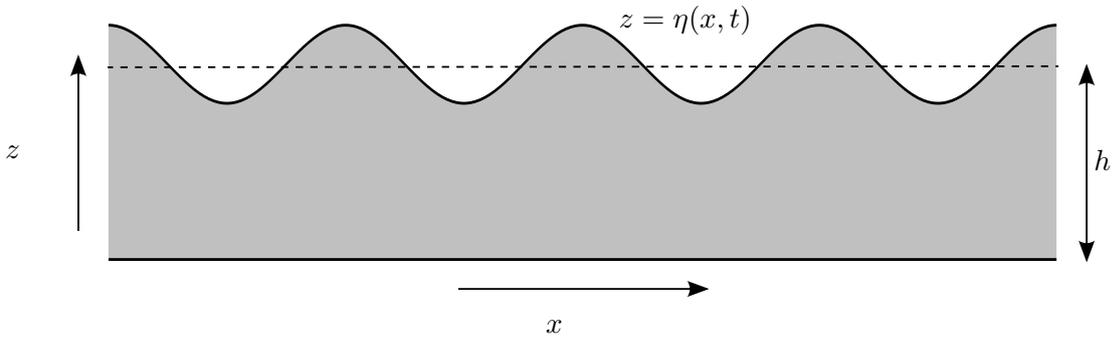}
\caption{\la{fig:ww} The domain for the water wave problem. Here $z=0$ is the equation of the surface for flat water, $z=-h$ is the flat bottom.
    }
\end{figure}
%

The main goal of the water wave problem is to understand the dynamics of the free surface $\eta(x,t)$. Thus it is convenient to recast the problem so as to involve only surface variables. Zakharov \cite{zakharov1} showed that the water wave problem is Hamiltonian with canonical variables $\eta(x,t)$ and $\varphi(x,t)=\phi(x,\eta(x,t),t)$. In other words $\varphi(x,t)$ is the velocity potential evaluated at the surface. Following \cite{craigsulem}, the Hamiltonian is written as

\beq
H=\frac{1}{2}\int_0^{2\pi}\left(\varphi G(\eta) \varphi+g \eta^2\right)dx,
\eeq

\no where $G(\eta)$ is the Dirichlet~$\ra$~Neumann operator: $G(\eta)\varphi=(1+\eta_x^2)^{1/2}\phi_n$, at $z=\eta(x,t)$. Here $\phi_n$ is the normal derivative of $\phi$. Using the water wave problem, $G(\eta)\varphi=\phi_z-\eta_x \phi_x=\eta_t$, which is the first of Hamilton's equations. The water wave problem for $\eta(x,t)$ and $\varphi(x,t)$ is

\beq\la{wweqs}
\eta_t=\dd{H}{\varphi}, ~~\varphi_t=-\dd{H}{\eta}.
\eeq

\begin{enumerate}

\item {\bf Quadratic Hamiltonian.} Since for our purposes, the linearization of (\ref{wweqs}-b) in a moving frame is required, it suffices to evaluate the Dirichlet~$\ra$~Neumann operator $G(\eta)$ at the flat surface $\eta=0$, resulting in

    \beq
    G(0)=-i \p_x \tanh(-ih \p_x).
    \eeq

The quadratic Hamiltonian $H_c^0$ is given by

\beq
H_c^0=c \int_0^{2\pi} \varphi \eta_x dx+\frac{1}{2}\int_0^{2\pi}\left(
\varphi (-i \tanh(-ih\p_x) \varphi_x+g \eta^2
\right) dx,
\eeq

\no giving rise to the linearized equations in a frame moving with velocity $c$:

\alpheqn

\begin{align}\la{linww}
\eta_t&=\dd{H_c^0}{\varphi}=c\eta_x-i \tanh(-ih \p_x)\varphi_x,\\
\varphi_t&=-\dd{H_c^0}{\eta}=c\varphi_x+g\eta.
\end{align}

\resetalpheqn

\item {\bf Dispersion relation.} The well-known dispersion relation \cite{vandenboek} for the water wave problem is immediately recovered from the linearized system (\ref{linww}-b) with $c=0$ (no moving frame), resulting in

\beq\la{wwdisp}
\omega^2=gk \tanh(kh).
\eeq

\no Note that the right-hand side of this expression is always positive. Thus there are two branches to the dispersion relation:

\beq\la{wwbranches}
\omega_{1,2}=\pm \sgn(k)\sqrt{gk\tanh(kh)}.
\eeq

\no Thus $\omega_1$ ($\omega_2$) corresponds to positive (negative) phase speed, independent of the sign of $k$.

\item {\bf Bifurcation branches.} Branches originate from $(c,\mbox{amplitude})=(\omega_1(1)/1,0)$ and $(c,\mbox{amplitude})=(\omega_2(1)/1,0)$. Without loss of generality, we focus on the first branch for which the phase speed $\sqrt{g\tanh(h)}$ is positive. This allows for a straightforward comparison of our results with those for the Whitham equation, in Example~\ref{ex:whitham}.

\item {\bf Stability spectrum.} The elements of the spectrum are given by

\alpheqn
\begin{align}\nonumber
\lambda_{n,1}^{(\mu)}&=-i\,\Omega_1(n+\mu)\\\la{wwspectrum}
&=i(n+\mu)\sqrt{g\tanh(h)}-i\,\sgn(n+\mu)\sqrt{g(n+\mu)\tanh(h(n+\mu))},\\
\nonumber
\lambda_{n,2}^{(\mu)}&=-i\,\Omega_2(n+\mu)\\
&=i(n+\mu)\sqrt{g\tanh(h)}+i\,\sgn(n+\mu)\sqrt{g(n+\mu)\tanh(h(n+\mu))}.
\end{align}
\resetalpheqn

\no The $\sgn(n+\mu)$'s may be omitted in these expressions, as the same set of spectral elements is obtained.

\item {\bf Collision condition.} The condition \rf{cancollision} is easily written out explicitly, but for our purposes it suffices to plot $\Omega_1(k+n)$ and $\Omega_2(k+n)$, for different values of $n$. This is done in Fig.~\ref{fig:wwcollision}b with $g=1$ and $h=1$. Although only the first collision is visible in the figure (all intersection points are horizontal integer shift of each other and correspond to the same value of $\mu$ and $\lambda_{n,j}^{(\mu)}$), it is clear from the curves shown that many collisions occur. The figure is qualitatively the same for all finite values of depth.

%
\begin{figure}[tb]
\begin{tabular}{cc}
\def\svgwidth{3.2in}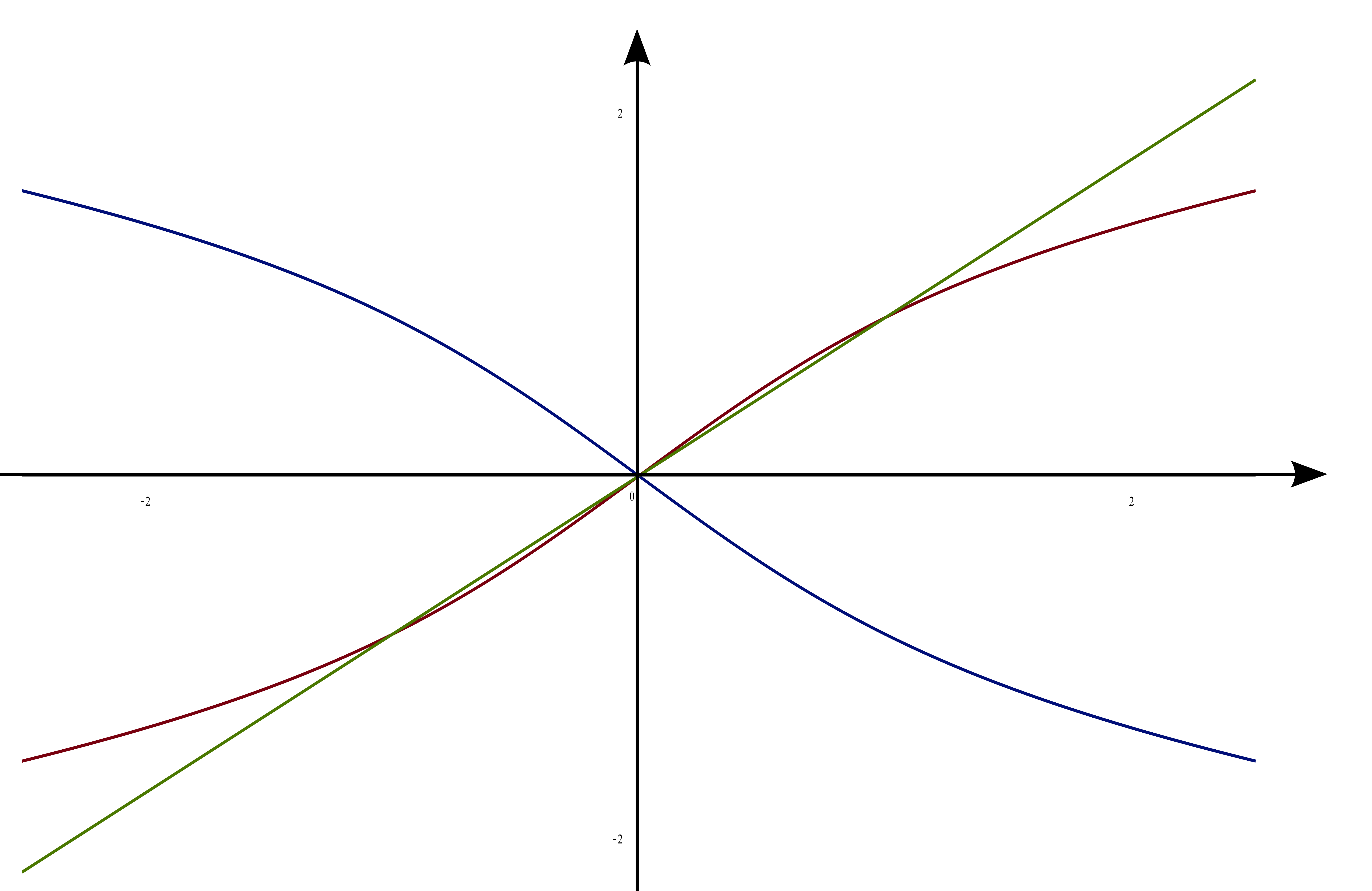 &
\def\svgwidth{2.8in}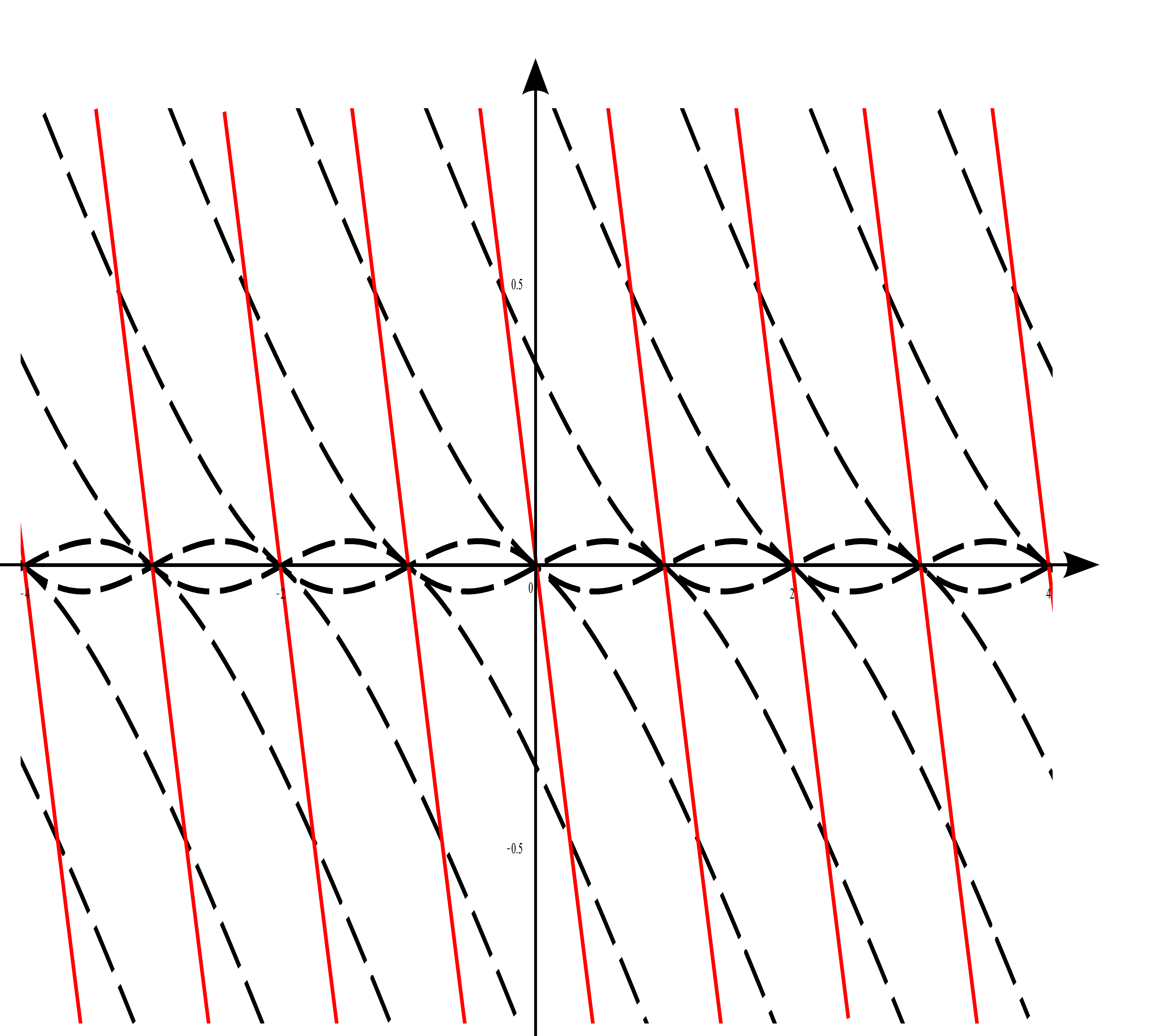\\
(a) & (b)\\
\multicolumn{2}{c}{
\def\svgwidth{3.2in}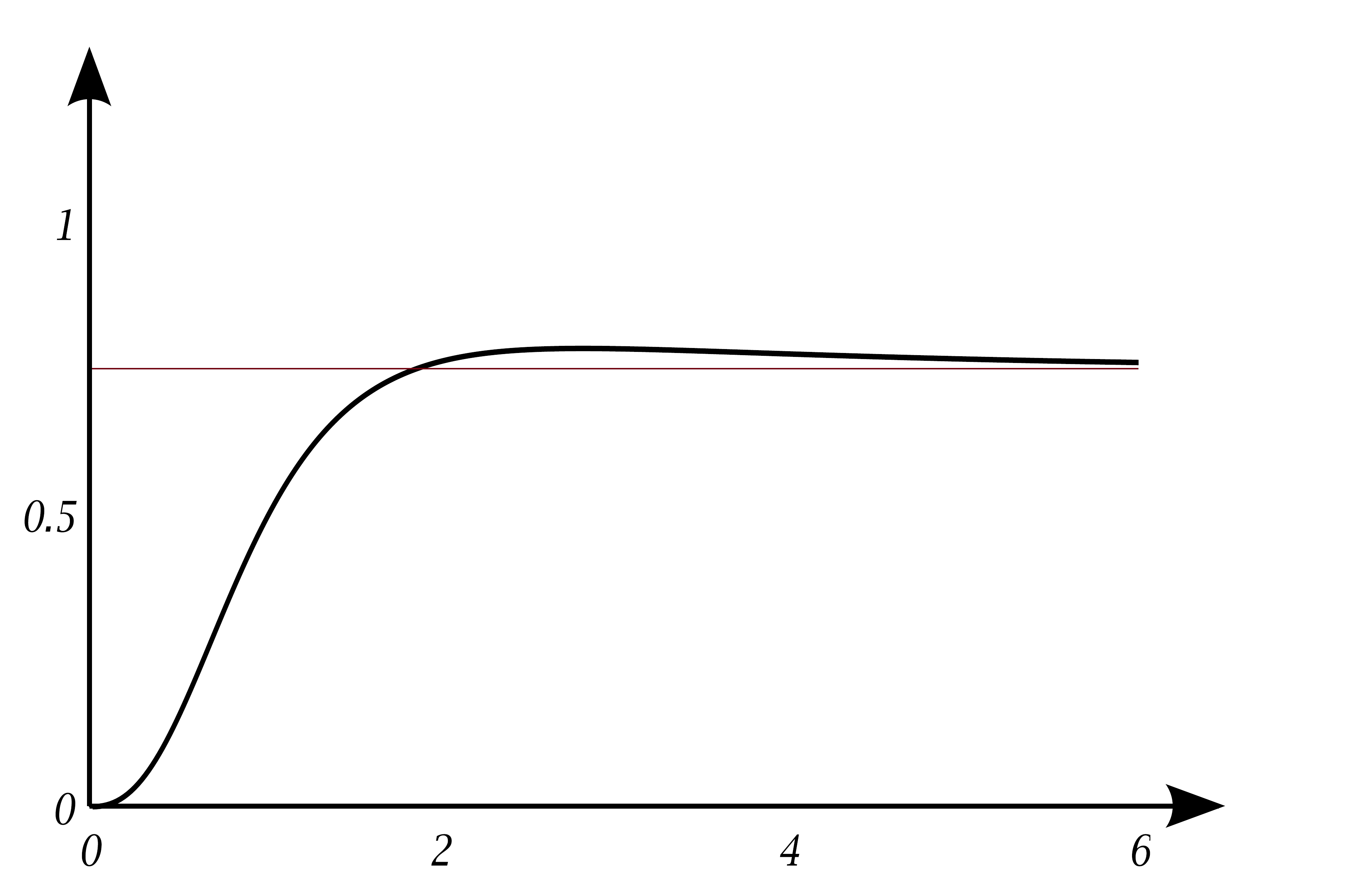
}\\
\multicolumn{2}{c}{(c)}
\end{tabular}
\caption{\la{fig:wwcollision} (a) The two branches of the dispersion relation for the water wave problem ($g=1$, $h=1$). The line through the origin has slope $\omega_1(1)/1$, representing the right-hand side of \rf{cancollision}. (b) The two families of curves $\Omega_1(k+n)$ (red, solid) and $\Omega_2(k+n)$ (black, dashed), for various (integer) values of $n$, illustrating that many collisions occur away from the origin. (c) The origin of the high-frequency instability closest to the origin as a function of depth $h$.
    }
\end{figure}
%

\item {\bf Krein signature.} The conditions \rf{sym1} and \rf{sym2} become

\beq\la{sym1ww}
\omega_{l_1}(n_1+\mu)\omega_{l_2}(n_2+\mu)g^2<0,
\eeq

\no and

\beq\la{sym2ww}
\omega_{l_1}(n_1+\mu)\omega_{l_2}(n_2+\mu)\sum_{{j_1}=1}^\infty \alpha_{j_1-1}h^{2j_1-1} (n_1+\mu)^{2j_1}\sum_{{j_2}=1}^\infty \alpha_{j_2-1}h^{2j_2-1} (n_2+\mu)^{2j_2}<0,
\eeq

\no respectively. Here the coefficients $\alpha_j$ are related to the Bernoulli numbers \cite{dlmf}, as they are defined by the Taylor series

\beq\la{tanh}
\tanh(z)=\sum_{j=0}^\infty \alpha_j z^{2j+1}, ~~~|z|<\pi/2.
\eeq

\no Because of the finite radius of convergence of this series, \rf{sym2ww} is only valid for small values of the wave numbers $n_1+\mu$ and $n_2+\mu$, but it is possible to phrase all results in terms of $\tanh$ directly, avoiding this difficulty. For instance, using \rf{tanh}, \rf{sym2ww} may be rewritten as

\begin{align*}
&&\omega_{l_1}(n_1+\mu)\omega_{l_2}(n_2+\mu)\frac{\omega_{\alpha}^2(n_1+\mu)}{g}
\frac{\omega_{\beta}^2(n_2+\mu)}{g}&<0\\
&\Rightarrow& \omega_{l_1}(n_1+\mu)\omega_{l_2}(n_2+\mu)&<0,
\end{align*}

\no in agreement with \rf{sym1ww}. The indices $\alpha$, $\beta\in \{0,1\}$ are irrelevant since $\omega_\alpha$ and $\omega_\beta$ both appear squared. This serves to illustrate that for specific examples one of the two criteria \rf{cankrein1} and \rf{cankrein2} (or \rf{sym1} and \rf{sym2} for even systems) may be significantly easier to evaluate, although they are equivalent.

Thus all collision points are potential origins of high-frequency instabilities. It appears from the numerical results in \cite{deconinckoliveras1} that the bubble of non-imaginary eigenvalues closest to the origin contains the high-frequency eigenvalues with the largest real part. Thus for waves in shallow water $kh<1.363$ (no Benjamin-Feir instability) \cite{benjaminperiodic, whitham1, zakharovostrovsky}, these are the dominant instabilities. For waves in deep water ($kh>1.363$) the Benjamin-Feir instability typically dominates, although there is a range of depth in deep water where the high-frequency instabilities have a larger growth rate, see \cite{deconinckoliveras1}. The dependence on depth $h$ of the location on the imaginary axis from which the high-frequency bubble closest to the origin bifurcates is shown in Fig.~\ref{fig:wwcollision}(c), with $g=1$. As $h\rightarrow \infty $, the imaginary part of $\lambda\rightarrow 3/4$. This asymptote is drawn in Fig.~\ref{fig:wwcollision}(c) for reference. This figure demonstrates that for all positive values of the depth $h$, the instabilities considered are not modulational as they do not bifurcate away from the origin as the amplitude increases.

It was remarked in Example~\ref{ex:whitham} that no collisions are possible due to the concavity of the dispersion relation. As a consequence, all collisions away from the origin observed in Fig.~\ref{fig:wwcollision}b involve both branches of the dispersion relation, {\em i.e.}, they involve a solid curve and a dashed curve. This is easily seen from Fig.~\ref{fig:wwcollision}a: a parallel cord with abscissae of the endpoints that are integers apart is easily found by sliding a parallel cord away from the cord $((0,0),(1,\omega_1(1)))$ until the integer condition is met. This implies that $\omega_{l_1}(n_1+\mu)$ and $\omega_{l_2}(n_2+\mu)$ in the collision condition \rf{cancollision} have opposite sign and \rf{sym1ww} is always satisfied. Thus colliding eigenvalues of zero-amplitude water wave solutions {\em always} have opposite Krein signature. As a consequence, the necessary condition for the presence of high-frequency instabilities is met. In fact, it was observed in \cite{deconinckoliveras1} that all colliding eigenvalues give rise to bubbles of instabilities as the amplitude is increased.

Our general framework easily recovers the results of MacKay \& Saffman \cite{mackaysaffman}. There the set-up is for arbitrary amplitudes of the traveling wave solutions, but the results are only truly practical for the zero-amplitude case.

\end{enumerate}

{\bf Remark.} It follows from these considerations that the high-frequency instabilities present in the water wave problem are a consequence of counter-propagating waves as no such instabilities are present in the Whitham equation \rf{whitham}. Although it is often stated that the value of the Whitham equation lies in that it has the same dispersion relation as the water wave problem (see for instance \cite{whitham}), this is in fact not the case as it contains only one branch of the dispersion relation. Thus the equation does not allow for the interaction of counter-propagating modes, and as such misses out on much of the important dynamics of the Euler equations.

\section{A Boussinesq-Whitham equation}

The goal of this section is the introduction of a model equation that has the same dispersion relation as the Euler equations (\ref{eulera}-d) at the level of heuristics that led Whitham to the model equation \rf{whitham}. In other words, we propose a bidirectional Whitham equation, so as to capture both branches of the water wave dispersion relation. We refer to this equation as the Boussinesq-Whitham (BW) equation. It is given by

\beq\la{bw}
q_{tt}=N(q)+\p_x^2\left(\alpha q^2+\int_{-\infty}^\infty K(x-y)q(y)dy
\right),
\eeq

\no where

\beq
K(x)=\frac{1}{2\pi}\int_{-\infty}^\infty c^2(k)e^{ikx}dk,
\eeq

\no and $c^2(k)=g\tanh(kh)/k$, $\alpha>0$ for the water wave problem without surface tension. In \rf{bw}, $N(u)$ denotes the nonlinear terms, which are ignored in the remainder of this section. Since our methods focus on the analysis of zero amplitude solutions, the sign of $\alpha$ is not relevant in what follows. This equation is one of many that may stake its claim to the name ``Boussinesq-Whitham equation". Equation \rf{bw} is a ``Whithamized'' version of the standard Bad Boussinesq equation and it may be anticipated that it captures at least the small-amplitude instabilities of the water wave problem in shallow water. It should be remarked that the Bad Boussinesq equation is ill posed as an initial-value problem \cite{mckean1}, but it might be anticipated that the inclusion of the entire water-wave dispersion relation overcomes the unbounded growth that is present due to the polynomial truncation. We return to this at the end of this section.

Before applying our method to examine the potential presence of high-frequency instabilities of small-amplitude solutions of the BW equation, we need to present its Hamiltonian structure. Further, since \rf{bw} is defined as an equation on the whole line, a periodic analogue is required, as in Section~\ref{sec:whitham}.

It is easily verified that \rf{bw} is Hamiltonian with (non-canonical) Poisson operator \cite{mckean1}

\beq\la{bwj}
J=\left(
\ba{cc}
0 & \p_x\\
\p_x & 0
\ea
\right),
\eeq

\no and Hamiltonian

\beq\la{bwh}
H=\int_{-\infty}^{\infty} \left(\frac{1}{2} p^2+\frac{\alpha}{3} q^3 \right)dx+\frac{1}{2}\int_{-\infty}^{\infty} dx \int_{-\infty}^{\infty} dy\, K(x-y) q(x) q(y).
\eeq

\no Indeed, \rf{bw} can be rewritten in the form \rf{ham} with $u=(q,p)^T$.

To define a periodic version of \rf{bw}, let

\beq
K(x)=\frac{1}{L}\sum_{j=-\infty}^\infty c^2(k_j) e^{i k_j x},
\eeq

\no where $k_j=2\pi j/L$, $j\in \mathbb{Z}$. The periodic BW equation is obtained from \rf{ham}, using \rf{bwj} and \rf{bwh}, but with all $\pm$ infinities in the integration bounds replaced by $\pm L/2$, respectively. Since \rf{bw} has a Poisson operator \rf{bwj} that is different from those used in Sections~\ref{sec:scalar} and~\ref{sec:vector}, minor modifications to the use of the method are necessary.

\begin{enumerate}

\item {\bf Quadratic Hamiltonian.} Ignoring the contributions of the nonlinear term, the quadratic Hamiltonian in a frame of reference moving with speed $V$ is given by

\beq\la{bwhc}
H_V^0=\int_{0}^{2\pi}\left(V q p+\frac{1}{2}p^2 \right)dx+\frac{1}{2}\int_{0}^{2\pi} dx \int_{0}^{2\pi} dy\, K(x-y) q(x) q(y),
\eeq

\no where we have fixed the period of the solutions to be $L=2\pi$. The inclusion of the first term in \rf{bwhc} is one place where the effect of the different form for $J$ is felt, as its functional form is a direct consequence of the form of \rf{bwj}.

\item {\bf Dispersion Relation.} A direct calculation confirms that

\beq
\omega^2=gk \tanh(kh),
\eeq

\no which is, by construction, identical to the dispersion relation for the full water wave problem \rf{wwdisp}. This gives rise to two branches of the dispersion relation \rf{wwbranches}, corresponding to right- and left-going waves.

\item {\bf Bifurcation Branches.} Bifurcation branches for $2\pi$-periodic solutions start at $(V_{1,2},0)$, where the phase speeds $V_{1,2}$ are given by $V_{1,2}=\pm \sqrt{g\tanh(h)}$.

\item {\bf Stability Spectrum.} The stability spectrum elements are, again by construction, identical to those for the water wave problem, given in (\ref{wwspectrum}-b).

\item {\bf Collision Condition.} Given that the spectral elements are identical to those for the water wave problem, the collision condition is identical too. It is displayed in Fig.~\ref{fig:wwcollision}(a-b). Thus, collisions away from the origin occur. It remains to be seen whether these can result in the birth of high-frequency instabilities.

\item {\bf Krein Signature.} As for the canonical case of Section~\ref{sec:vector}, we use \rf{kreinhessian}. Thus we calculate the Hessian ${\cal L}_c$ of the Hamiltonian $H_c^0$.

    Let

    \beq
    c^2(k)=\sum_{j=0}^{\infty} \gamma_j k^{2j},
    \eeq

\no where $\gamma_j=gh^{2j+1}a_j$, with the coefficients $a_j$ defined in \rf{tanh}. A direct calculation gives that the Hamiltonian \rf{bwhc} is rewritten as

\beq
H_V^0=\frac{1}{2}\int_0^{2\pi} \left(p^2+Vqp+\sum_{j=0}^\infty \gamma_j q_{jx}^2\right)dx.
\eeq

\no Using this form of the Hamiltonian, the calculation of the Hessian is straightforward, leading to

\beq
{\cal L}_V=
\left(
\ba{cc}
\ds \sum_{j=0}^\infty \gamma_j (-1)^j \p_x^{2j} & V\\
V & 1
\ea
\right).
\eeq

Next, we compute the eigenvectors $v=(q,p)^T$. We have

\beq
\left(
\ba{c}
q\\p
\ea
\right)=e^{i \lambda t}
\left(
\ba{c}
Q(x)\\P(x)
\ea
\right),
\eeq

\no where $(Q,P)^T$ satisfies

\beq\la{bwevp}
\lambda \left(
\ba{c}
Q\\P
\ea
\right)= \left(
\ba{cc}
0 & \p_x\\
\p_x & 0
\ea
\right) {\cal L}_V \left(
\ba{c}
Q\\P
\ea
\right).
\eeq

This is a second place where the Poisson operator $J$ plays a crucial role as it affects the form of $v=(q,p)^T$ and thus the expression for the signature. One easily verifies that

\beq
\left(
\ba{c}
Q\\P
\ea
\right)=e^{i(n+\mu)x}\left(
\ba{c}
i(n+\mu)\\\lambda-i(n+\mu)V
\ea
\right)
\eeq

\no satisfies \rf{bwevp}.

We need to evaluate the sign of

\begin{align}\nonumber
&\left(
\ba{c}
Q\\P
\ea
\right)^\dagger
{\cal L}_V \left(
\ba{c}
Q\\P
\ea
\right)\\\nonumber
&=e^{-i(n+\mu)x}\left(
\ba{c}
-i(n+\mu) \\ -\lambda+i(n+\mu)V
\ea
\right)^T
\left(
\ba{cc}
\ds \sum_{j=0}^\infty \gamma_j (-1)^j \p_x^{2j} & V\\
V & 1
\ea
\right)\left(
\ba{c}
i(n+\mu)\\\lambda-i(n+\mu)V
\ea
\right)e^{i(n+\mu)x}\\\nonumber
&=\left(
\ba{c}
-i(n+\mu) \\ -\lambda+i(n+\mu)V
\ea
\right)^T
\left(
\ba{cc}
\ds \sum_{j=0}^\infty \gamma_j (n+\mu)^{2j} & V\\
V & 1
\ea
\right)\left(
\ba{c}
i(n+\mu)\\\lambda-i(n+\mu)V
\ea
\right)\\\nonumber
&=\left(
\ba{c}
-i(n+\mu)\\ i\omega(n+\mu)
\ea
\right)^T\left(
\ba{cc}
c^2(n+\mu) & V\\
V & 1
\ea
\right)\left(
\ba{c}
i(n+\mu)\\-i\omega(n+\mu)
\ea
\right)\\
&=2\omega\left(\omega-(n+\mu)V\right).
\end{align}

Let the signature associated with the first eigenvalue be the sign of $2\omega_{j_1}(\omega_{j_1}-(n_1+\mu)V)$, where $\omega_{j_1}$ is a function of $n_1+\mu$. Similarly, for the second eigenvalue, the signature is the sign of $2\omega_{j_2}(\omega_{j_2}-(n_2+\mu)V)$. Using the collision condition $\lambda_{n_1,j_1}^{(\mu)}=\lambda_{n_2,j_2}^{(\mu)}$, the product of these two expressions is $4\omega_{j_1}\omega_{j_2}(\omega_{j_2}-(n_2+\mu)V)^2$,
which is less than zero since collisions can only occur for eigenvalues associated with opposite branches of the dispersion relation, see Fig.~\ref{fig:wwcollision}b. It follows that, as in the water wave case, the signatures of colliding eigenvalues are always opposite, and the necessary condition for spectral instability is met. Thus, unlike the Whitham equation \rf{whitham}, the BW model \rf{bw} does not exclude the presence of high-frequency instabilities of small-amplitude solutions.

\end{enumerate}

The results obtained from the Krein signature calculations are confirmed by numerical results, see Fig.~\ref{fig:bwstuff}. Panel~(a) shows a numerically computed traveling wave solution of the BW equation \rf{bw}. This solution is computed using a cosine collocation method with 60 points, as for the Whitham equation, see Fig.~\ref{fig:wstuff} \cite{SanfordKodamaCarterKalisch}. For the solution plotted, $c\approx 1.049815$. The second panel in the first row displays the spectrum computed using Hill's method with 100 modes and 20000 values of the Floquet parameter, using an interpolation of the solution profile. This panel shows the presence of a large number of apparent instabilities, most with small growth rate, in the neighborhood of the imaginary axis. The third panel shows a zoom of the region around the origin, revealing a modulational instability. This is expected, since such an instability is also present for the Whitham equation, see Section~\ref{ex:whitham}. The fourth panel zooms in on the first bubble of instabilities centered on the positive imaginary axis, revealing a shape and location that is consistent with the Krein collision theory presented here.

%
\begin{figure}[hbt]
\begin{tabular}{cc}
\hspace*{-1.0in}\def\svgwidth{3.2in}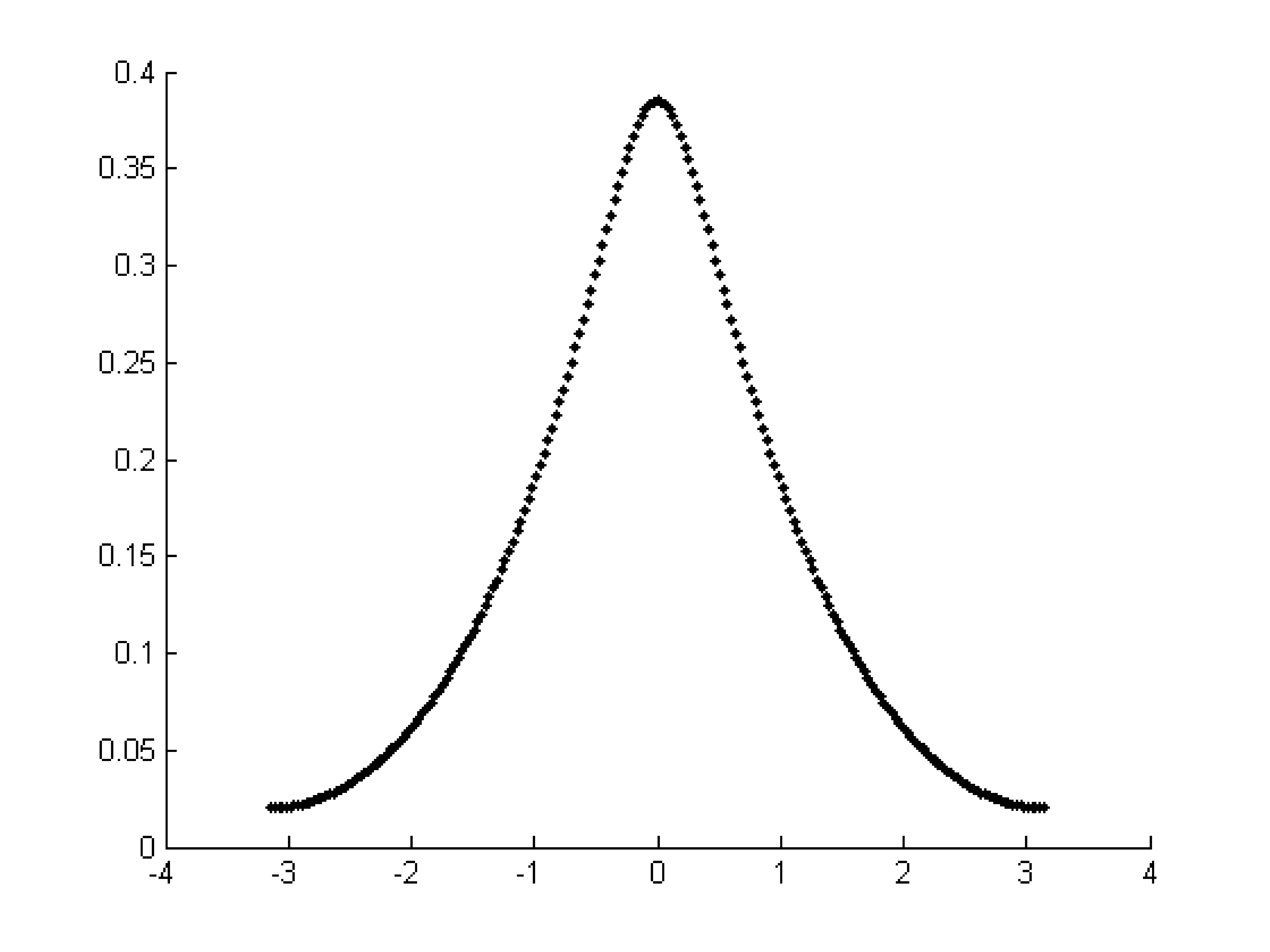 &
\hspace*{-1.14in}\def\svgwidth{4in}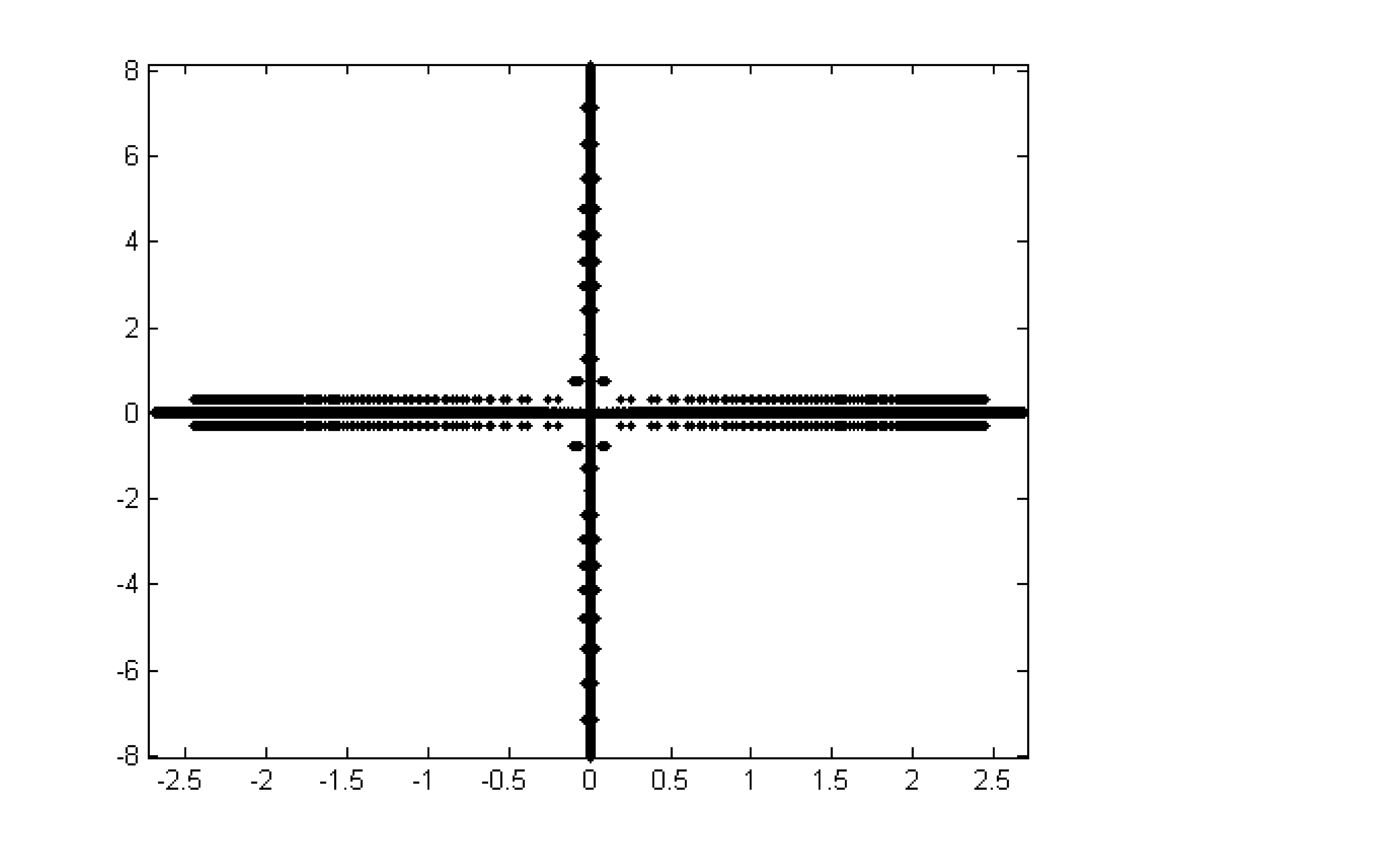\\
\hspace*{-0.7in}(a) & \hspace*{-1.5in}(b)\\
\hspace*{-0.3in}\def\svgwidth{4in}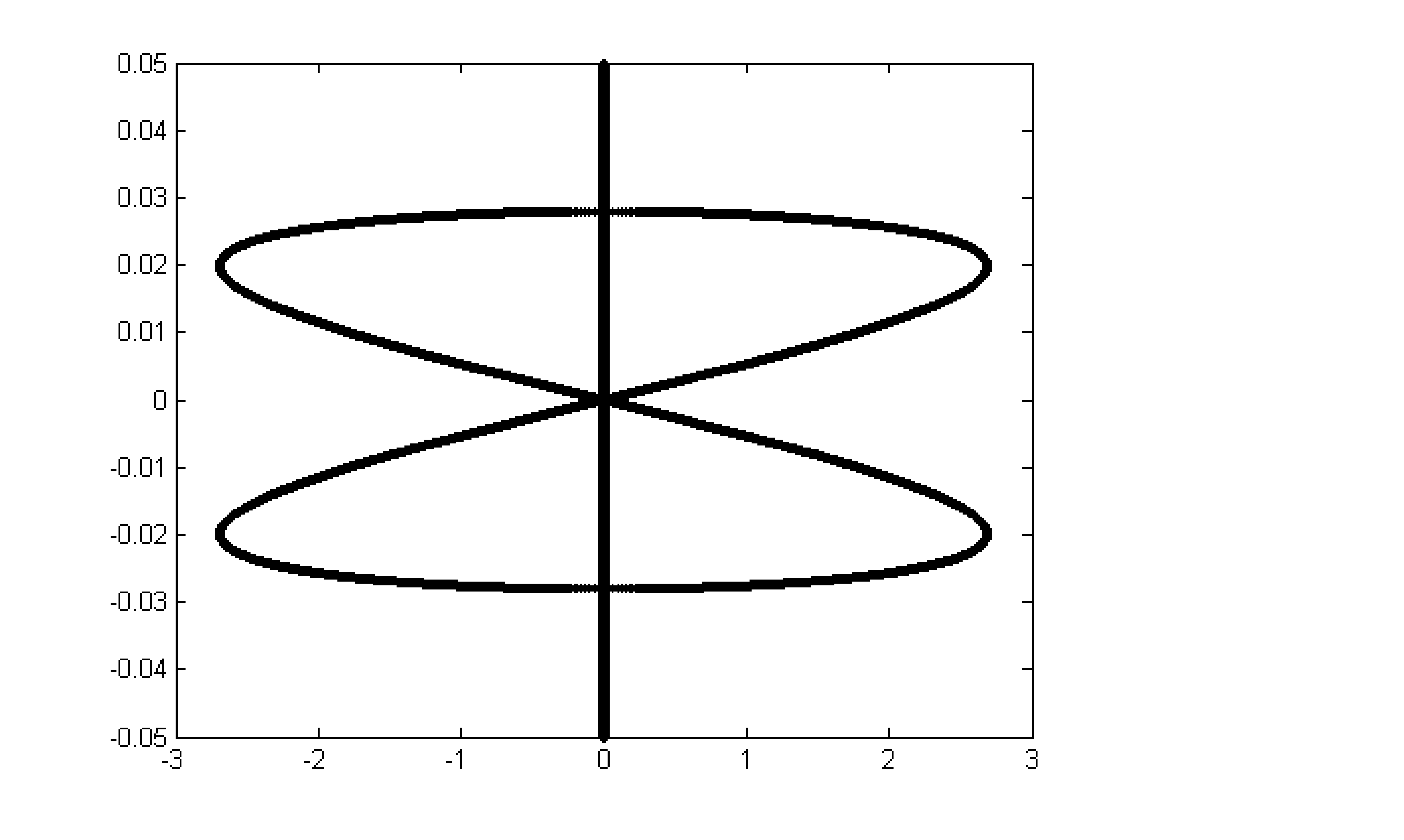 &
\hspace*{-1.14in}\def\svgwidth{4in}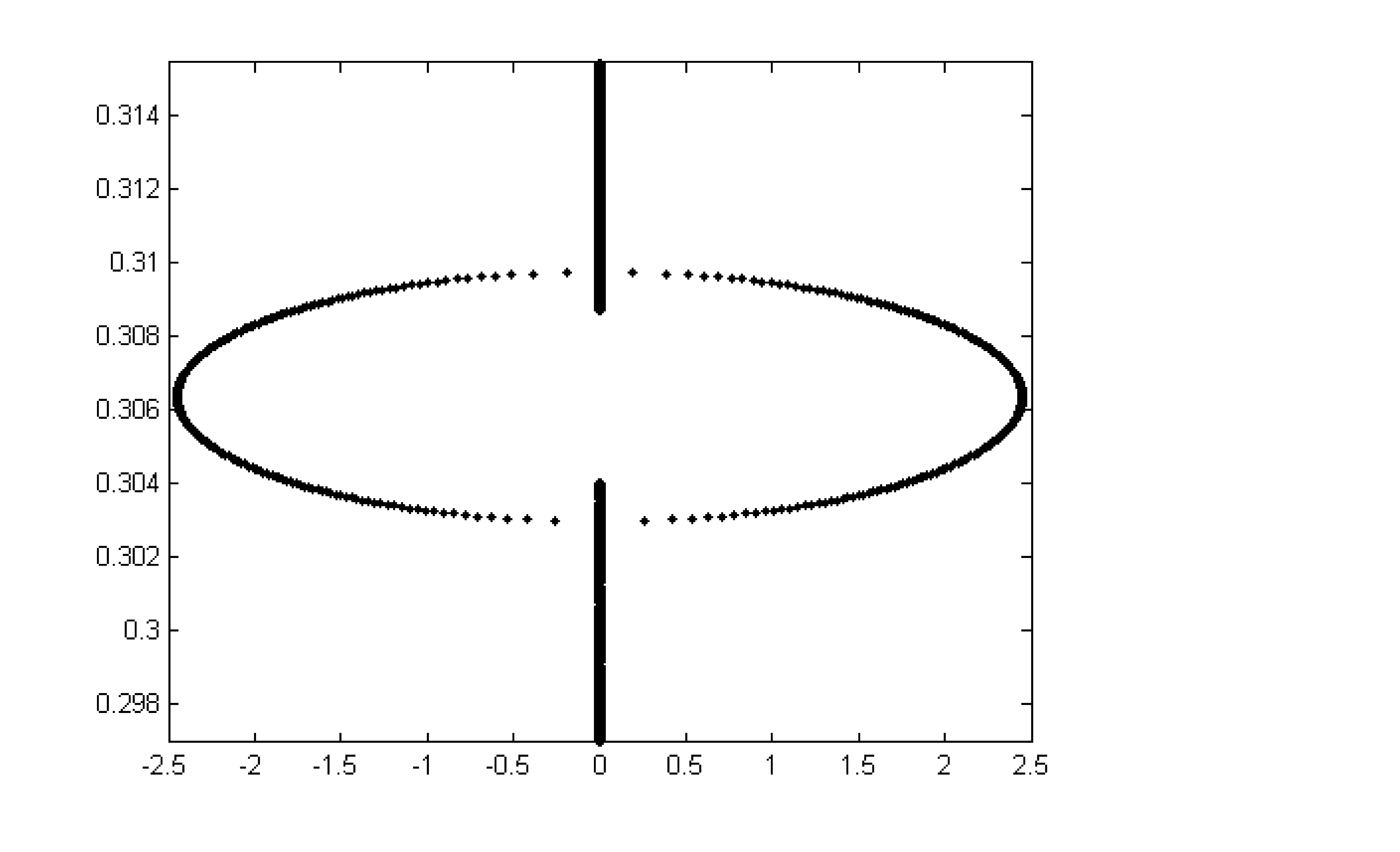\\
\hspace*{-0.7in}(c) & \hspace*{-1.5in}(d)\\
\end{tabular}
\caption{\la{fig:bwstuff} (a) A small-amplitude traveling wave solution of the Boussines-Whitham equation \rf{bw} with $c\approx 1.0498515$. (b) The numerically computer stability spectrum. (c) A blow-up of the stability spectrum in a neighborhood of the origin. (d) A blow-up of the stability spectrum around what appears as a horizontal segment visible in (b) immediately above the longest segment appearing horizontal. More detail is given in the main text.
    }
\end{figure}
%

\vs

The wave form displayed in Fig.~\ref{fig:bwstuff}(a) does not have zero average, unlike the one shown in Fig.~\ref{fig:wstuff}, for reasons explained here. Let us examine the stability of a flat-water state $q=a$ (constant), $p=0$. Thinking of the BW equation as an approximation to the water wave problem where the flat-water state is neutrally stable (spectrum on the imaginary axis), independent of the reference level of the water, the neutral stability of this state is desired in the context of the BW system as well. However, the BW system is easily checked to not be Galilean invariant, thus the average value of the solution may be important.

Linearizing the system around the flat water state $(q,p)=(a,0)$ results in a linear system with constant coefficients, whose dispersion relation is given by $\omega^2=2ak^2+k^2c^2(k)$. This results in two branches for the dispersion relation: $\omega_{1,2}=\pm k \sqrt{c^2(k)+2a}$. It follows that if $a>0$ then both $\omega_1$ and $\omega_2$ are real, resulting in neutral stability, since the stability eigenvalue and the frequency $\omega(k)$ are different by a factor of $i$. On the other hand, if $a<0$, it follows that both $\omega_1$ and $\omega_2$ are imaginary for sufficiently large $k$, since $\lim_{|k|\ra \infty}c(k)=0$. This leads to the dynamics of the flat-water state with $a<0$ to not only be unstable, but to be ill-posed, as the growth rate of the instability $\ra \infty$ as $|k|\ra \infty$. Thus Whitham-izing the Bad Boussinesq equation and incorporating the full water wave dispersion relation does not remove the illposedness of the problem. Rather it alters it where negative constant solutions experience unbounded growth, unlike positive constant solutions.

It is observed numerically that this behavior of perturbed constant solutions is carried over to nonconstant solutions: solutions of negative average display the same illposed behavior described above, with stability spectra that have unbounded real part. In contrast, the spectra of solutions of positive average have bounded real part, as in Fig.~\ref{fig:bwstuff}. Annoyingly, the illposed behavior extends to numerical solutions constructed to have zero average. Presumably this is a consequence of numerical error, as higher accuracy numerical experiments display narrower spectra whose real part tends to infinity more slowly.

We may summarize our findings on the BW equation as follows. The equation was constructed as a bi-directional Whitham equation so as to truly have the same linear dispersion relation as the water wave problem. Even though the BW has a different Poisson structure than the water wave problem, we find that periodic solutions of the BW are susceptible to high-frequency instabilities originating from Krein collisions at the exact same locations on the imaginary axis as the water wave problem. On the other hand, we have not attempted to quantify whether the resulting growth rates are comparable to those for the water wave problem. Further, the illposedness of the equation for solutions of negative average is a significant strike against its potential use in applications. Nevertheless, it appears possible to design more equations like the BW equation, possessing the exact same dispersion relation as the water wave problem and with it all its high-frequency instabilities, without the equation dynamics being illposed for any important class of solutions.

\section*{Acknowledgements}

We wish to thank Richard Kollar for interesting discussions. John Carter helped our understanding of the Whitham equation and Mat Johnson was a part of our initial investigations on the Boussinesq-Whitham equation.  This work was supported by the National
Science Foundation through grant NSF-DMS-1008001 (BD) and in part by the EPSRC under grant EP/J019569/1 and NSERC (OT). Any opinions, findings,
and conclusions or recommendations expressed in this material are
those of the authors and do not necessarily reflect the views of the
funding sources.

\bibliographystyle{plain}

\bibliography{mybib}

\end{document}